\documentclass[11pt]{article}
\usepackage{graphicx,psfrag}
\usepackage{amssymb,amsfonts,amsthm}
\usepackage{amsmath,amsthm,amssymb}
\usepackage{amsfonts}
\usepackage{relsize}
\usepackage[T1]{fontenc}
\usepackage[utf8]{inputenc}
\usepackage{graphicx}

\DeclareMathOperator*{\bigast}{\raisebox{-0.6ex}{\scalebox{2.5}{$\ast$}}}

\newcommand{\mb}{\mathbf}
\newcommand{\geod}{\mathrm{geod}}
\newcommand{\la}{\langle}
\newcommand{\ra}{\rangle}
\newcommand{\SL}{\mathrm{SL}}
\newcommand{\MCG}{\mathcal{MCG}(S_{g,p})}
\newcommand{\M}{\mathcal{M}}
\newcommand{\cc}{\mathcal{C}}
\renewcommand{\P}{\mathcal{P}}
\newcommand{\F}{\mathcal{F}}
\newcommand{\D}{\mathcal{D}}
\newcommand{\si}{\!\!\sim}
\newcommand{\co}{\mathfrak{c}}
\newcommand{\rc}{\mathfrak{rc}}
\newcommand{\calr}{{\mathcal R}}
\newcommand{\calc}{{\mathcal C}}
\newtheorem{theorem}{Theorem}[section]
\newtheorem{lemma}[theorem]{Lemma}

\newtheorem{cor}[theorem]{Corollary}
\newtheorem{proposition}[theorem]{Proposition}
\theoremstyle{definition}
\newtheorem{definition}[theorem]{Definition}

\newtheorem{remark}[theorem]{Remark}
\newtheorem{prob}[theorem]{Problem}

\newcommand{\call}{{\mathcal L}}

\newcommand{\dist}{{\mathrm{dist}}}
\newcommand{\ddd}{{\mathcal D}}

\newcommand{\nn}{{\mathcal N}}

\newcommand {\N}{\mathbb{N}} 
\newcommand {\Z}{\mathbb{Z}}            
\newcommand {\R}{\mathbb{R}} 
\newcommand {\C}{\mathbb{C}} 
\newcommand {\supp}{\mathrm{supp}}

\newcommand {\me}{\medskip}

\newcommand {\iv}{^{-1}}

\newcommand{\prop}{\mathrm{prop}}
\begin{document}

\title{The Rapid Decay property and centroids in groups}

\author{Mark Sapir\thanks{The
research was supported in part by the NSF grants DMS 1418506 and
DMS	1318716, and a BSF grant.}}

\date{}
\maketitle

\me

\begin{abstract}
This is a survey of methods of proving or disproving the Rapid Decay property in groups. We present a centroid property of group actions on metric spaces. That property is a generalized (and corrected)
version of the  ``(**)-relative hyperbolicity" from \cite{DS} and implies the Rapid Decay (RD) property. We show that several properties which are known to imply RD also imply the centroid property. Thus uniform lattices in many semi-simple Lie groups, graph products of groups, Artin groups of large type and the mapping class groups have the (relative) centroid property. We also present a simple ``non-amenability-like" property that follows from RD, and give an easy example of a group without RD and without any amenable subgroup with superpolynomial growth.
\end{abstract}

\tableofcontents

\section{Introduction}\label{s:0}

$\gimel(n)$

Recall that a length function on a group $G$ is a map $L$ from $G$
to the set of non-negative real numbers $\R_+$ satisfying:
\begin{itemize}
    \item[(1)] $L(gh)\leq L(g)+ L(h)$  for all $g,h\in G$ ;
   \item[(2)] $L(g)=L(g\iv)$  for all $g\in G$ ;
    \item[(3)] $L(1)=0$;
\end{itemize}
A length function is called \emph{proper} if
\begin{itemize}
    \item[(4)] For every $r>0$ the set $\{g\in G\mid L(g)\le r\}$ is finite.
\end{itemize}

For a finitely generated group any length function $L$ is \emph{dominated} by any word length function $L_S$ induced by a finite generating set $S$, that is, $L(g)\le CL_S(g)$ for some constant $C$ and every element $g\in G$. Indeed, if $w=s_1s_2\ldots s_n$, where $s_i\in S$, then $L(w)\le L(s_1)+L(s_2)+\ldots +L(s_n)\le Cn$ where $C=\max(L(s_1),\ldots, L(s_n))$.

If $G$ acts on a metric space $X$ by isometries, $x_0\in X$, then the function $L_{x_0}\colon g\to \dist(x_0,g\cdot x_0)$ is a length function (easy to check).

Let $G$ be a countable group equipped with a length function $L$ and the corresponding pseudo-distance $\dist(a,b)=L(a\iv b)$ (this would be a distance function if $L(g)=0\to g=1$).
The analytic definition of property RD introduced by Haagerup and Jolissaint (see Jolissaint's paper \cite{J} or Valette's book \cite{V}) is the following.
For every $s\in \R$ {\it the
Sobolev space of order $s$ with respect to $L$} is the set $H_L^s
(G)$ of functions $\phi$ on $G$ such that the function $(1+L)^s\phi$
is in $l^2 (G)$. The space of {\it{rapidly decreasing functions on
$G$ with respect to $L$}} is the set $H_L^\infty (G)=\bigcap_{s\in
\R }H_L^s (G)$.

The {\it{group algebra of $G$}} over $\C$, denoted by $\C G$, is the set of
functions with finite support on $G$.

With every element $g\in G$ we can associate the linear {\em
convolution operator} $\phi\mapsto g*\phi$ on $l^2(G)$, where
$$g*\phi(h)=\phi(g\iv h).$$ This is just the left regular representation of $G$ on $l^2(G)$, it can be extended to a representation of $\mathbb{C}G$ on
$l^2(G)$ by linearity. This representation is faithful and every
convolution operator induced by an element of $\C G$ is bounded.
Therefore we can identify $\C G$ with a subspace in the space of
bounded operators $\mb{B}(l^2(G))$ on $l^2(G)$. For every $x\in
\C G$ we denote by $\|x\|_*$ its operator norm, that is
$$
\| x\|_* = \sup \{ \|x*\phi\| \; ;\; \|\phi\| =1\}\, .
$$


\begin{definition}
The group $G$ is said to have {\it{the RD property with respect to the length-function $L$}} if the inclusion of $\C G$ into the reduced $C^*$-algebra $C_r^*(G)$ of $G$ extends to a continuous inclusion of $H_L^\infty
(G)$ into $C_r^*(G)$.
\end{definition}

One can reformulate the property RD in the following way involving only real valued non-negative functions with finite supports (see \cite{DS}).
\begin{definition}\label{d:rd}
Let $\phi$ be a function $G\to \R_+$ with finite support $\supp(\phi)=\{g\in G, \phi(g)\ne 0\}$. The ($l_2$-)norm is defined as usual: $||\phi||=\sqrt{\sum_{g\in G} \phi(g)^2}$. The maximal length of an element from the support of $\phi$ will be denoted by $\prop(\phi)$ and is called the \emph{propagation} of $\phi$. If $\phi, \psi$ are two functions with finite supports, then $\phi*\psi$ is the function $G\to \R_+$ defined by $\phi*\psi(k)=\sum_{g\in G}\phi(g)\psi(g\iv k)$. We say that $G$ has the property RD if there is a polynomial \footnote{Here and below all coefficients of all polynomials are assumed to be non-negative, so that the polynomials are strictly increasing on $\R_+$.} $P(x)$ such that for every positive number $r$ and every functions $\phi, \psi\colon G\to \R_+$ with finite supports such that $\prop(\phi)\le r$, and we have

$$||\phi*\psi||^2\le P(r)||\phi||^2||\psi||^2.$$
\end{definition}

Note \cite{DS} that if a group satisfies property RD with respect to some length function $L$, then it satisfies  RD with respect to any length function that dominates $L$. In particular, a finitely generated group satisfies property RD if an only if it satisfies RD with respect to the word length function (induced by a finite generating set).

Property RD turned out to be important in several areas of mathematics, from analytic K-theory to C*-algebras to random walks on Cayley graphs of groups. Most notably, groups having property RD ``very often" satisfy the Baum-Connes conjecture (without coefficients), hence the Novikov conjecture, etc. \cite{V1,V}.

Many classes of groups are known to satisfy RD. After Haagerup proved it for the free groups \cite{Ha}, Jolissaint and de la Harpe proved it for all Gromov hyperbolic groups \cite{J,dlH}. Non-uniform lattices in higher rank semi-simple Lie groups do not have RD.  One of the most stimulating conjectures in the area is the conjecture of Valette \cite[Conjecture 7]{V} that all uniform lattices in semi-simple Lie groups should have property RD. That conjecture is still wide open even for $\SL_4(\R)$ in spite of a lot of efforts.  By results of Ramagge, Robertson, Steger \cite{RRS}, Lafforgue \cite{L}, Chatterji \cite{Ch} and Talbi \cite{T} we know that every uniform lattice in $\SL_3(K)$ where $K$ is a field $\R$ or $\mathbb{C}$ or a ring of quaternions or octonions, and in many  direct products of such Lie groups and Lie groups of rank 1 have property RD. More recently Chatterhi, Ruane \cite{CR} and Dru\c tu and myself \cite{DS} proved property RD for groups that are relatively hyperbolic with respect to groups with RD, Behrstock and Minsky \cite{BM} proved it for the mapping class groups of surfaces, Ciobanu, Holt and Rees proved RD for large type Artin groups and for graph products of groups with RD \cite{CHR1}, \cite{CHR}.

Even though the classes of groups which are known to have RD are quite different, the methods of proofs are ``asymptotically similar". The reason why the free groups have RD is that every geodesic triangle on a tree has a center which belongs to every side of the triangle. For Cayley graphs of hyperbolic groups, a center of a triangle may not belong to all three sides, but it is at bounded distance from all three sides (this is Rips' definition of hyperbolic groups). For triangles in the Cayley graphs of relatively hyperbolic groups \cite{DS} and in symmetric spaces of Lie groups such as  $\SL_3(\R)$ \cite{L} every triangle has an ``inscribed" nice and relatively small triangle from certain family of triangles (properties (*) and $(K_\delta)$ below). Chatterji and Ruane used clouds of centers \cite{CR}, and Ciobanu, Holt and Rees \cite{CHR} used a condition which can be interpreted as a center-like condition (see below). The goal of this paper is to present an easy to formulate and check ``centroid" condition which follows from the centroid-type conditions used before and implies property RD. Such an attempt was made at the end of our paper \cite{DS}. There we formulated our property ``(**)-relative hyperbolicity". But that property is not general enough and  the definition of (**) in \cite{DS} contains errors.

The centroid and relative centroid properties introduced in Section \ref{s:2} below can be considered as the ``true (**)-relative hyperbolicity". We present proofs that both the centroid property and the relative centroid property with respect to sets of triples satisfying RD imply  property RD. Although the proofs are similar and the first result follows from the second one, we present both proofs for educational reasons (the first proof is much easier and more ``natural"). Also algebraists like me do not often get the pleasure of using the Cauchy-Schwarz inequality in their papers, and it is used twice in each proof. In Section \ref{s:3} we shall show that the (relative) centroid property follows from several ``center-like" properties studied earlier. In Section \ref{s:4}, we present a combinatorial consequence  of property RD, and give an example of a group without RD and without amenable subgroups of superpolynomial growth (the question of existence of such groups was discussed at the AIM workshop on property RD (Palo Alto, January 23 to January 27, 2006).

{\bf Acknowlegement.} I would like to thank Jason Behrstock, Indira Chatterji, Laura Ciobanu, Cornelia Dru\c tu, Paul Jolissaint, Mitchel Kleban, Bogdan Nica, Denis Osin and Sarah Rees for helpful conversations. I would also like to thank Yair Minsky for pointing out a mistake in \cite{DS}.

\section{The centroid and relative centroid properties}\label{s:2}

\subsection{The centroid property and RD}

\begin{definition}\label{d:0} Let $G$ be a countable group acting by isometries on a metric space $(X,\dist)$, $x_0\in X$ with point stabilizers finite of uniformly bounded sizes. We assume that $L$  is the length function defined by $L(g)=\dist(x_0,g\cdot x_0)$ (as in Section \ref{s:0}).  Let $\co$ be a map from the set of pairs $G^2=G\times G$ to $X$. We can view $G$ as embedded into $X$ (by the map $g\to g\cdot x_0$), a pair $(g,k)\in G\times G$ as the vertices of triangle $(x_0,g\cdot x_0,k\cdot x_0)$, and $c=\co(g,k)$ as a ``center" of that triangle.  We say that $G$ and $\co$ satisfy the \emph{centroid property} if for some polynomial $P(x)$ we have

\begin{enumerate}
\item[($c_1$)] For every $k\in G$ and every $r>0$ the number of elements in the set $\{\co(g,k), L(g)\le r\}$ does not exceed $P(r)$.
\item[$(c_2)$] For every $g\in G$ the number of elements in the set $\{\co(g,k), k\in G\}$ does not exceed $P(L(g))$.
\item[$(c_3)$] For every $h\in G$ the number of elements in the set $\{g\iv\co(g,gh), L(g)\le r\}$ does not exceed $P(r)$.
\end{enumerate}
In this case $X$ will be called the \emph{space of centroids} of $G$  and $\co$ will be called the \emph{centroid map}.
\end{definition}

It is obvious that every group satisfying the centroid property with respect to a length function $L$ also satisfies this property with respect to any length function that dominates $L$.

\begin{remark} Note that if $\co$ can be equivariantly extended from triangles $(1,g,k)$ to arbitrary triangles $(a,b,c)\in G^3$ and, as a map on $G^3$, $\co(a,b,c)$ is invariant under the permutations of variables, then $(c_3)$ follows from $(c_1)$. Indeed, $\co(g,gh)=\co(1,g,gh)=g\co(g\iv, 1,h)=g\co(1,g\iv,h)=g\co(g\iv, h)$.
\end{remark}

\begin{theorem}\label{t1} The centroid property implies property RD.
\end{theorem}

\proof  Let $\phi, \psi\colon G\to \R_+$ be two functions with finite supports and the propagation of $\phi$ is equal to $r$. We need to estimate $||\phi*\psi||^2$ from above. By definition $$||\phi*\psi||^2=\sum_k\left(\sum_{g\in\supp(\phi)} \phi(g)\psi(g\iv k)\right)^2.$$ We can rewrite this sum as

\begin{equation}\label{1}
\sum_k\left(\sum\limits_{\substack{c \\ \exists g\in\supp(\phi)\colon\\ c=\co(g,k)}}\sum\limits_{\substack{g \\ c=\co(g,k)}}\phi(g)\psi(g\iv k)\right)^2.\end{equation}


 Let us use the following corollary of the  Cauchy-Schwarz inequality which is true for all real $a_i\in \R$:

\begin{equation}\label{5}\left(\sum_{i=1}^n a_i\right)^2=\left(\sum_{i=1}^n 1\cdot a_i\right)^2\le \sum_{i=1}^n 1^2\sum_{i=1}^n a_i^2 \le n\sum a_i^2.\end{equation}

We can apply (\ref{5}) to the first inner sum in (\ref{1}). By ($c_1$), for each $k\in G$ the number of possible points $c$ such that $c=\co(g,k)$, where $L(g)\le r$, does not exceed  $P(r)$.  Therefore the sum in (\ref{1}) does not exceed

\begin{equation}\label{2}\begin{array}{l}P(r)\sum\limits_k\sum\limits_{\substack{c \\ \exists g\colon c=\co(g,k)}}\left(\sum\limits_{\substack{g \\ c=\co(g,k)}}\phi(g)\psi(g\iv k)\right)^2. \\ \hskip 1 in
\end{array}
\end{equation}

Applying the Cauchy-Schwarz inequality to the inner sum in (\ref{2}), we deduce that (\ref{2}) does not exceed

\begin{equation}\label{8}P(r)\sum\limits_{k}\sum\limits_{\substack{c}}\left(\sum\limits_{\substack{g\colon\\  c=\co(g,k)}}\phi(g)^2\right)\left(\sum\limits_{\substack{h\colon\\ c=\co(kh\iv,h),\\ L(g)\le r}}\psi(h)^2\right)\end{equation}
(we denoted $g\iv k$ by $h$).

Let us compute the number of times the expression $\phi(g)^2\psi(h)^2$ for given $g, h\in G$ occurs in the expansion of (\ref{8}). It is easy to see that it is equal to the number of pairs $(k,c)$ such that  for some $g_1,h_1\in G$, $L(g_1)\le r$, we have

\begin{equation}\label{7} c=\co(g,gh_1)=\co(g_1, g_1h), g_1h=gh_1=k.\end{equation}

Let us fix $g,h$, $L(g)\le r$. Then the number of possible points $c=\co(g,gh_1)$ is at most $P(r)$ by ($c_2$). Now, in addition to fixing $g,h$, let  us fix $c=\co(g,gh_1)=\co(g_1,g_1h)$ for some $g_1,h_1$ where $g_1h=gh_1$ and estimate the number of possible elements $k=g_1h$. We have
$k\iv \cdot c=h\iv\cdot (g_1\iv\cdot \co(g_1,g_1h))$. By ($c_3)$ the number of points $g_1\iv\cdot \co(g_1,g_1h)$ with $L(g_1)\le r$ is at most $P(r)$. Since orders of the point stabilizers of the action of $G$ on $X$ are uniformly bounded by some constant $K$, we have that the
number of possible elements $k$ (given $g,h,c$) does not exceed $K P(r)$. Therefore the number of pairs $(k,c)$ for any given $g, h$ does not exceed $KP^2(r)$. Hence
$$||\phi*\psi||^2\le KP(r)^3||\phi||^2||\psi||^2$$ and property RD follows.
\endproof

\subsection{The relative centroid property}

The centroid property is a generalization (and correction) of the (**)-relative hyperbolicity with respect to the trivial subgroup from \cite{DS}. The full (**)-relative hyperbolicity can be generalized too in a very similar manner.

First we need to slightly generalize the property RD (this generalization can be found in \cite{L}).

We say that an action of a group $G$ on a metric space is almost free if the diameters of all point stabilizers are uniformly bounded from above.

Let $X$ be a metric space and $G$ be a group acting on $X$ almost freely.
For every $\mb{x},\mb{y}\in X^2/G$, the product $\mb{x}\mb{y}$ is defined as $\mb{x}\mb{y}=\{(\alpha,\gamma), \exists \beta: (\alpha,\beta)\in \mb{x}, (\beta,\gamma)\in \mb{y}\}$.

In general $\mb{x}\mb{y}$ is a union of orbits of $G$ but if the action is free, the situation is better.

\begin{lemma}\label{l56}
If the action of $G$ in $X$ is free, then for every $\mb{x},\mb{y}\in X^2/G$ we have $\mb{x}\mb{y}\in X^2/G$.
\end{lemma}
\proof
Indeed, if $(\alpha,\gamma), (\alpha',\gamma')\in \mb{x}\mb{y}$, then for some $\beta,\beta'\in X$, we have $(\alpha,\beta)\in \mb{x},(\beta,\gamma)\in \mb{y}, (\alpha',\beta')\in \mb{x}, (\beta',\gamma')\in \mb{y}$. Then there exists $g,g'\in G$ such that $g\cdot (\alpha,\beta)=(\alpha',\beta'), g'\cdot (\beta,\gamma)=(\beta',\gamma')$. Then $g\cdot \beta=g'\cdot \beta$. Hence $g=g'$ since the action is free, and $g\cdot(\alpha,\gamma)=(\alpha',\gamma')$. Thus $\mb{x}\mb{y}$ coincides with an orbit of $G$ in $X^2$.
\endproof

A $G$-orbit from $X^2$ containing a pair $(\alpha,\beta)$ will be denoted by $[\alpha,\beta]$.

\begin{lemma}\label{l57} If the action of $G$ on $X$ is free, then every $k\in G$
is uniquely determined by two points $\alpha,\beta\in X$ and the orbit $[\alpha,k\cdot \beta]$.
\end{lemma}

\proof Indeed, since the action of $G$ on $X$ is free, for every $\alpha\in X$, the orbit contains at most one pair of the
form $(\alpha,\delta)$.
Thus, given $\alpha$ and $[\alpha,k\cdot\beta]$, we can determine $k\cdot \beta$ uniquely. Since $\beta$ is given, we can determine $k$.
\endproof

%

A function $\phi$ from $X^2\to \R_+$ is called $G$-invariant if $\phi(g\cdot x,g\cdot y)=\phi(x,y)$ for every $g\in G, x,y\in X$. In that case $\phi$ induces a function $X^2/G\to \R_+$ which we shall denote by $\phi$ also. We say that  $G$-invariant function $\phi$ has \emph{finite support} if
its support consists of finite number of $G$-orbits $\mb{x}\in X/G$.
Another fact that we will be using is the following.

Given two functions $\phi,\psi\ X^2\to \R_+$ and a subset $T\subset X^3$, we can define the convolution of $\phi, \psi$ relative to $T$:

$$\phi*_T\psi(x,y)=\sum_{z\in X, (x,z,y)\in T} \phi(x,z)\psi(z,y).$$

Note that if $\phi, \psi, T$ are $G$-invariant, then $\phi*_T\psi$ is $G$-invariant and the support of $\phi*_T\psi$ consists of some (possibly not all) orbits  $\mb{x}\mb{y}$ where $\mb{x}\in X^2/G$ is in the support of $\phi$, $\mb{y}\in X^2/G$ is in the support of $\psi$.

The \emph{propagation} $\prop(\phi)$ of a function $\phi\colon X^2\to \R_+$ is the maximal $\dist(x,y)$ for all pairs $(x,y)$ from the support of $\phi$. The \emph{norm} $||\phi||$ for a $G$-invariant function $\phi\colon X^2\to \R_+$ is the norm of the induced function on $X^2/G$ (that is the square root of the sum of squares of values of $\phi$ on the orbits of $G$). Thus

$$
||\phi||^2=\sum_{\mb{x}\in X^2/G} \phi(\mb{x})^2.
$$
So for every two functions $\phi,\psi\colon X^2\to \R$ with finite supports

\begin{equation}\label{56} ||\phi*\psi||^2=\sum\limits_{\mb{z}\in X^2/G}\left(\sum\limits_{\mb{x},\mb{y}\in X^2/G: \mb{z}\subseteq \mb{x}\mb{y}} \phi(\mb{x})\psi(\mb{y})\right)^2.
\end{equation}

We say that a $G$-invariant subset $T\subseteq X^3$ satisfies property RD if there is a polynomial $P(r)$ such that for every two $G$-equivariant functions $\phi,\psi\colon X^2\to \R_+$ with finite support and $\prop(\phi)\le r$ we have

$$||\phi*_T\psi||^2\le P(r) ||\phi||^2||\psi||^2.$$

\begin{remark}\label{rk:rd} Clearly a group $G$ with length function $L$ has property RD if and only if it has RD relative to the set $G^3$ of all triples (the centroid space is $G$ itself). Moreover if $G$ is a group with length function $L$ and $H$ is a subgroup of $G$ with induced length function $L_H$, then $H$ has RD with respect to the length function $L_H$ if and only if the set of triples $G\cdot H^3=\{(gh_1, gh_2, gh_3), g\in G, h_1, h_2, h_3\in H\}$ has RD.
\end{remark}

\begin{definition} Let $G$ be a group acting almost freely on a metric space $X$ (the \emph{space of centroids}). Let $T_1,\ldots,T_n$ be $G$-invariant subsets of $X^3$. We say that the group $G$ has relative centroid property with respect to $T_1,\ldots,T_n$ if there exists a function $\rc\colon G^2\to \bigcup_{i =1}^m T_i$ and a polynomial $P(r)$ such that

($rc_1$) For every $k\in G$ the number of pairs $(\alpha,\gamma)\in X^2$ such that for some $g\in G, L(g)\le r$ and $\beta\in X$, $\rc(g,k)=(\alpha,\beta,\gamma)$ does not exceed $P(r)$.

($rc_2$) For every $g\in G$ the number of pairs $(\alpha,\beta)\in X^2$ such that for some $k\in G$ and $\gamma\in X$,
$\rc(g,k)=(\alpha,\beta,\gamma)$ does not exceed $P(L(g))$.

($rc_3$) For every $h\in G$ the number of pairs $g\iv\cdot (\beta,\gamma)\in X^2$ such that for some $g\in G, L(g)\le r$ and $\alpha\in X$,
$\rc(g,gh)=(\alpha,\beta,\gamma)$ does not exceed $P(r)$.

($rc_4$) For every $g,k\in G$, if $(\alpha,\beta,\gamma)=\rc(g,k)$, then $$\dist(\alpha,\beta)\le P(L(g)), \dist(\alpha,\gamma)\le P(L(k)), \dist(\beta,\gamma)\le P(L(g\iv k)).$$
\end{definition}

Note that the
centroid property is the same as the relative centroid property with respect to the diagonal $T$ of $X^3$ where $X$ is the centroid space.

\begin{theorem}\label{t10} If $G$ as above has a relative centroid property with respect to sets $T_1,\ldots, T_m\in X^3$ which have property RD, and the length function is proper, then $G$ has property RD,
\end{theorem}

\proof The proof is similar to the proof of Theorem \ref{t1} (and to the proof of \cite[Theorem 3.1]{DS}). In the proof, we assume, to simplify formulas, that $m=1$. The case of $m>1$ is very similar and is left to the reader. Let us denote $T_1$ by $T$.

First assume that the action of $G$ on $X$ is free. We shall deal with almost free actions later.

For every triple $x=(\alpha,\beta,\gamma)$ and every subset $S$ of $\{1,2,3\}$ let $\pi_S(x)$ be the projection of $x$ onto the coordinates from $S$. For example, $\pi_{13}(x)=(\alpha,\gamma)$.

To simplify formulas, we need the following notation.

\me


For every $k\in G$, let $\ddd_k$ be the set of triples $(\alpha,[\alpha,\gamma],k\iv \cdot\gamma)\in X\times X^2/G\times X$
such that for some $g\in G,\gamma\in X$, $(\alpha,\beta,\gamma)=\rc(g,k)$.
For every $d\in \ddd_k$ let $\calc_d$ be the set of triples $([\alpha,\beta],[\beta,\gamma],g\iv\cdot\beta)\in X^2/G\times X^2/G\times X$ such that  $(\alpha,[\alpha,\gamma],k\iv\cdot \gamma)=d$ and $\rc(g,k)=(\alpha,\beta,\gamma)$.
%

The sets $\pi_{13}(\ddd_k)$, $\pi_1(\ddd_k)$, $\pi_3(\ddd_k)$ are denoted by $\call\calr_k$, $\call_k$ and $\calr_k$ respectively.

We denote by $\ddd$ the union of all $\ddd_k, k\in G$.
The sets $\call\calr$, $\call$, $\calr$ are defined similarly. By Lemma \ref{l57}, the element $k$ is uniquely determined by any triple $(\alpha, [\alpha,\gamma], k\iv\cdot \gamma)$. Thus there is a natural map $\eta$ from the set of all triples $(\alpha, [\alpha,\gamma],k\iv\cdot \gamma)$ to $k$. It takes each $\ddd_k$ to $k$ (in particular the sets $\ddd_k$ for different $k$ are disjoint).
Note that if $(\alpha,\beta,\gamma)=\rc(g,k)$, then

\begin{equation}\label{19}\begin{array}{l}
\eta(\alpha,[\alpha,\beta],g\iv\cdot \beta)=g, \\
\eta(g\iv \beta,[\beta,\gamma],k\iv\cdot \gamma)=g\iv k,\\
 \eta(\alpha,[\alpha,\gamma],k\iv\cdot \gamma)=k.
\end{array}
\end{equation}

By Lemma \ref{l56} we also have

\begin{equation}\label{19.5}
[\alpha,\beta][\beta,\gamma]=[\alpha,\gamma].
\end{equation}

We can rewrite $||\phi*\psi||^2$ as follows:

\begin{equation}\label{10}\sum\limits_k\left(\sum\limits_{\substack{(\alpha,\gamma)\in X^2}}
\sum\limits_{\substack{g \\ (\alpha,\gamma)=\pi_{13}\rc(g,k)}}\phi(g)\psi(g\iv k)\right)^2.\end{equation}

We can apply (\ref{5}) to the first inner sum in (\ref{10}). Using ($rc_1$) we deduce that the
sum in (\ref{10}) does not exceed

$$
P(r)\sum\limits_k\sum\limits_{\substack{(\alpha,\gamma)}}\left(\sum\limits_{\substack{g \\ (\alpha,\gamma)=\pi_{13}\rc(g,k)}}\phi(g)\psi(g\iv k)\right)^2.
$$

By (\ref{19}) this does not exceed

\begin{equation}\label{20}P(r)\sum\limits_k\sum\limits_{d=(\alpha,
\mb{z},\kappa)\in \ddd_k}
\left(\sum\limits_{\substack{(\mb{x},\mb{y},\delta)\in \calc_d}}\phi(\eta(\alpha,\mb{x},\delta))\psi(\eta(\delta,\mb{y},\kappa))\right)^2
\end{equation}
(here $\kappa$ denotes $k\iv\cdot\gamma$, $\delta$ denotes $g\iv \cdot\beta$).

The inner sum in (\ref{20}) can be rewritten as

\begin{equation}\label{25}
\sum\limits_{(\mb{x},\mb{y})\in \pi_{12}(\calc_d)}\left(\sum\limits_{\delta: (\mb{x},\mb{y},\delta)\in \calc_d}\phi(\eta(\alpha,\mb{x},\delta))\psi(\eta(\delta,\mb{y},\kappa))\right)
\end{equation}


Applying the Cauchy-Schwarz inequality to the inner sum, we deduce that the number in (\ref{25}) does not exceed

\begin{equation}\label{27}
\sum\limits_{(\mb{x},\mb{y})\in \pi_{12}(\calc_d)}\left(\sum\limits_{\delta: (\mb{x},\mb{y},\delta)\in \calc_d}\phi(\eta(\alpha,\mb{x},\delta))^2\right)^\frac{1}{2}\left(\sum\limits_{\delta: (\mb{x},\mb{y},\delta)\in \calc_d}\psi(\eta(\delta,\mb{y},\kappa))^2\right)^\frac{1}{2}
\end{equation}

Fixing $k, \alpha, \kappa$ (and thus fixing $d$), we define two functions $\Phi, \Psi\colon X^2/G\to \R_+$ as follows. For every $\mb{x}, \mb{y}\in X^2/G$ we define
$\Phi(\mb{x})$ as
$$
\left(\sum\limits_{\delta: \exists \mb{y} (\mb{x},\mb{y},\delta)\in \calc_d}\phi(\eta(\alpha,\mb{x},\delta))^2\right)^\frac{1}{2}
$$
and $\Psi(\mb{y})$ as

$$\left(\sum\limits_{\delta: \exists \mb{x} (\mb{x},\mb{y},\delta)\in \calc_d}\psi(\eta(\delta,\mb{y},\kappa))^2\right)^\frac{1}{2}
$$

 Both functions are  with finite supports. Moreover $\prop(\Phi)\le P(r)$ by ($rc_4$). Then (\ref{27}) does not exceed

$$\sum\limits_{\mb{x},\mb{y}:\mb{x}\mb{y}=\mb{z}} \Phi(\mb{x})\Psi(\mb{y}).$$

Thus (\ref{25}) does not exceed $(\Phi*\Psi)(\mb{z})$ and
(\ref{20}) does not exceed

\begin{equation}
\label{30}
P(r)\sum\limits_k\sum\limits_{(\alpha,\kappa)\in\call\calr_k}||\Phi*_T\Psi||^2
\end{equation}

By property RD for $T$, there exists a polynomial $P'(r)$ such that (\ref{30}) does not exceed

$$P(r)P'(P(r))\sum\limits_k\sum\limits_{(\alpha,\kappa)\in\call\calr_k} ||\Phi||^2||\Psi||^2.$$
which can be rewritten as

$$P(r)P'(P(r))\sum\limits_k\sum\limits_{(\alpha,\kappa)\in \call\calr_k}\left(\sum\limits_{\mb{x},\delta}\phi(\eta(\alpha,\mb{x},\delta))^2\right)
\left(\sum\limits_{\mb{y},\delta}\phi(\eta(\delta,\mb{y},\kappa))^2\right)$$
which does not exceed

\begin{equation}\label{100}
P(r)P'(P(r))\sum\limits_{(\alpha,\kappa)\in \call\calr}\left(\sum\limits_{\mb{x},\delta}\phi(\eta(\alpha,\mb{x},\delta))^2\right)
\left(\sum\limits_{\mb{y},\delta}\phi(\eta(\delta,\mb{y},\kappa))^2\right).
\end{equation}
This without the factor $P(r)P'(P(r))$ can be estimated from above by
$$
\left(\sum\limits_{\substack{(\alpha,\mb{x},\delta)\\ \exists g,k,\gamma: (\alpha,g\cdot\delta,\gamma)=\rc(g,k),\\ \mb{x}=[\alpha,g\cdot\delta]}} \phi(\eta(\alpha,\mb{x},\delta))^2\right)
\left(\sum\limits_{\substack{(\delta,\mb{y},\kappa)\\ \exists k,g\alpha: (\alpha,g\cdot\delta,k\cdot\kappa)=\rc(g,k),\\ \mb{y}=[g\cdot\delta,k\cdot \kappa]}} \psi(\eta(\delta,\mb{y},\kappa))^2\right)
$$

By Lemma \ref{l57}, for every $g\in G$ the number of times $\phi(g)^2$ appears in the  first sum in this expression is at most the number of pairs $(\alpha,\beta)=\pi_{12}\rc(g,k)$ where $k$ runs over $G$. This number does not exceed $P(r)$ by $(rc_2)$. Similarly the number of times $\phi(h)^2$ appears in the second sum  is at most the number of pairs $g\iv \pi_{23}\rc(g,gh)$ as $g$ runs over the set of elements of $G$ of length at most $r$. This number is at most $P(r)$ by $(rc_3)$. Thus (\ref{100}) does not exceed $$P(r)^3P'(P(r))||\phi||^2||\psi||^2.$$



Hence
$$||\phi*\psi||^2\le P(r)^3P'(P(r))||\phi||^2||\psi^2||.$$

Now let us assume that the action of $G$ on $X$ is almost free. Let $\dist_X$ be the metric on $X$. For every $x\in X$ let $G_x$ be the stabilizer of $x$ in $G$. Choose a point $x_C$ in every orbit $C$ of $G$ in $X$.
We can define a metric $\dist_C$ to $G$ by $\dist_C(g,h)=\dist(g\cdot x_C,h\cdot x_C)+1$ provided $g\ne h$, and $\dist_C(g,h)=0$ otherwise.  
The metric space $G$ with this metric will be denoted by $G_C$. Now let us consider the disjoint union $Y=\sqcup_C G_C$. We define the metric $\dist_Y$ on $Y$ as follows. If two points $g,h$ are in the same $G_C$, then $\dist_Y(g,h)=\dist_C(g,h)$, if $g\in G_C, h\in G_{C'}$, $C\ne C'$, then set $\dist_Y(g,h)=\dist_X(g\cdot x_C,h\cdot x_{C'})+1$. The group $G$ acts on $Y$ by $g\cdot h=gh$. We leave it to the reader to check that the action is free and by isometries. \footnote{The construction of the space $Y$ presented in \cite[Page 335]{CR}, is not
complete because an action of $G$ on $Y$ is not defined there, and it is not at all clear how to define such an action. A better
explanation (similar to the one we give here) can be found in \cite[Remark 2.15]{Chat}.} There is a natural map $\rho$ from $Y$ to $X$ which takes each $s\in G_C$ to $s\cdot x_C$. This map is $G$-equivariant, for every $y,y'\in Y$ we have
$$\dist_X(\rho(y), \rho(y'))\le \dist_Y(y,y')\le \dist_X(\rho(y), \rho(y'))+1$$
and for every $x\in X$ the diameter of the set $\rho\iv(x)$ does not exceed 1 and hence the number of elements in $\rho\iv(x)$ does not exceed a uniform constant $K$.

Let $T'=\rho\iv(T)$. Then $T'$ is a $G$-invariant subset of $Y^3$. It is easy to check that $T'$ satisfies property RD.
Let $\rc$ be the relative centroid map $G^2\to X^3$. Then for every $g,k$ let $\rc'(g,k)$ be any triple from $\rho\iv(\rc(g,k))$.

It is very straightforward to check that $\rc'$ satisfies properties $(rc_1)$-$(rc_4)$. Since the action of $G$ on $Y$ is free, we can apply the result we have already proved.
\endproof

Recall \cite{J} that a group has RD if and only if its subgroup of finite index has RD and if and only if factor-group over a finite normal subgroup has RD. Here are the analogs of these results for the relative centroid property.

\begin{definition}\label{de:rd}
We say that $G$ has a relative centroid property with respect to subgroups $H_1,\ldots, H_m$ if it has the relative centroid property with centroid space $G^3$ (with the natural action by $G$) and the sets of triples $T_i=G\cdot H_i^3, i=1,\ldots,m$ (see Remark \ref{rk:rd}).
\end{definition}

\begin{cor} Suppose that $H$ is a finite index subgroup of $G$. Then $G$ has a relative centroid property with respect to $H$.
\end{cor}

\proof Indeed, let $\{x_1,\ldots, x_n\}$ be representatives of left cosets of $H$ in $G$ and $\{y_1,\ldots, y_n\}$ be representatives of the right cosets of $H$ in $G$. Then  for every pair $(g,k)\in G^2$, let $g=x_ih, g\iv k=h'y_j$ for some $h,h'\in H, 1\le i,j\le n$.  We set $\rc(g,k)=
(x_i, x_ih, x_i hh')\in G\cdot H^3$. It is a straightforward exercise that $\rc$ satisfies properties $(rc_1)-(rc_4)$.
\endproof

\begin{cor} Suppose that $H=G/N$ where $N$ is a finite normal subgroup of $G$. Suppose that $H$ has a relative centroid property with respect to a centroid space $X$ and sets of triples $T_1,\ldots,T_n$. Then $G$ has a relative centroid property with respect to $X$ and $T_1,\ldots,T_n$.
\end{cor}

\proof Indeed, the almost free action of $H$ on $X$ induces an almost free action of $G$ on $X$. The centroid map $\rc$ on $G$ is obtained as a composition of the centroid map from $H$ and the natural homomorphism $G\to H$. It is easy to see that $\rc$ satisfies the conditions $(rc_1)-(rc_4)$.
\endproof

\section{Examples of groups with the centroid property}\label{s:3}

\subsection{The mapping class group of an oriented surface}

The following theorem is essentially proved by Behrstock and Minsky in \cite{BM}. Nevertheless the formulations of Theorems 1.2 and 3.2 of \cite{BM} contain mistakes (and the proof of Theorem 3.2 contains a mistake too)\footnote{\label{page}The mistake in \cite[Theorem 1.2]{BM} is that an equivariant map $\kappa\colon \MCG^3\to\MCG$ (\cite[Condition 2 of Theorem 1.2]{BM}) cannot be invariant under all permutations of the arguments (\cite[Condition 1 of Theorem 1.2]{BM}). Indeed, if $g$ is an element of order 3 (there are elements of order 3 in $\MCG$), then $\kappa(1,g,g^2)=g\kappa(g^2,1,g)=g\kappa(1,g,g^2)$, a contradiction. The same mistake was made in \cite[Section 4]{DS} where we gave an informal definition of property (**) (as was pointed out to us by Yair Minsky). It did not affect the main results of \cite{DS}.}.

\begin{theorem} \label{th:mcg} The  mapping class group $\MCG$ of an orientable surface of genus $g$ and $p$ punctures where $3g-3 +p \ge 1$ satisfies the centroid property (and hence has property RD).
\end{theorem}

\proof We shall use the notation and terminology from \cite{BM}. Only the centroids in the sense of \cite{BM} we shall call \emph{BM-centroids} to avoid confusion with the centroids from the centroid property.  Let $\M$ be the marking graph of $S_{g,p}$. Then $M$ is a locally finite graph and $\MCG$ acts on $\M$ by isometries, properly and co-compactly. Therefore the stabilizers of points of $\M$ in $\MCG$ are finite and their sizes are uniformly bounded, so the action is almost free.

We shall view $\M$ as the centroid space for $\MCG$.
Fix a marking $\mu\in \M$ with trivial stabilizer in $\MCG$. Consider the length function on $G$ defined by $L(g)=\dist(\mu,g\cdot \mu)$. Take any pair of elements $g,k\in \MCG$.
Then let $\co(g,k)$ be a BM-centroid of the triple $(\mu,g\cdot \mu,k\cdot\mu)$ (see \cite[The first paragraph of the proof of Theorem 3.2]{BM}).

Let us prove that $\co$ satisfies the conditions ($c_1$), ($c_2$), ($c_3$). First note that $\co(g,k)$ is at uniformly bounded distance from a geodesic $[\mu,g\cdot\mu]$ and also from some geodesics $[g\cdot\mu, k\cdot\mu]$ and $[\mu, k\cdot\mu]$. Therefore $\co(g,k)$ is contained in the intersection of the $\Sigma$-hulls of pairs of points $(\mu, g\cdot\mu)$, $(g\cdot\mu, k\cdot\mu)$, $(\mu, k\cdot\mu)$ (see \cite[The second paragraph of Section 4.1]{BM}).

($c_1$) This follows from \cite[Part 4 of Theorem 1.2]{BM} whose proof is not affected by the errors mentioned above.


($c_2$) immediately follows from \cite[Theorem 4.2]{BM}.

($c_3$) Let $r>0$ and $g,h\in \MCG$, $L(g)=r$. By construction $c=\co(g,gh)$ is the BM-centroid of the triangle $(\mu,g\cdot \mu,gh\cdot\mu)$. Also by construction, $g\iv\cdot c$ is a BM-centroid of the triangle $(g\iv\cdot \mu, \mu, h\cdot\mu)$ which is the same as the triangle $(\mu, g\iv \cdot \mu, h\cdot\mu)$. Since the set of possible centroids of any  triangle in $\M$ has uniformly bounded diameter (by \cite[Theorem 2.9]{BM}), the number of possible points of the form $g\iv\cdot \co(g,gh)$ (as $g$ varies) does not exceed $k_3r^\xi$ for some uniform constant $k_3$ where $\xi=3g-3+p$.
\endproof

%
%
%
%
%

\subsection{The Chatterji-Ruane property}

\begin{definition} [See \cite{CR}] Let $G$ be a group acting freely\footnote{This condition can be replaced by ``almost freely" as shown in \cite{Chat}, see also the end of the proof of Theorem \ref{t10}.} by isometries on a metric space $(X,\dist)$.
such that there is a $G$-equivariant map
$C \colon X\times X \to \mathcal{P}(X)$ (where $\mathcal{P}(X)$ is the set of all subsets of $X$),
$(x, y)\to  C(x, y)$,
satisfying the following (for any
$x, y, z \in X$).
\begin{itemize}
\item[(i)] $x\in C(x,y)$.
\item[(ii)] $C(x, y) \cap C(y, z) \cap C(z, x) \ne
\emptyset$.
\item[(iii)] There is a polynomial $R$ such that for any $r\in \R_+$, the
cardinality of $C(x, y) \cap  B(x, r)$ is bounded above by $R(r)$ where $B(x,r)$ is the ball of radius $r$ and center $x$ in $X$.
\item[(iv)] There is a polynomial $Q$ such that if $\dist(x, y) \le r$, then the diameter
of $C(x, y)$ is bounded by $Q(r)$.
\end{itemize}
Then we shall say that $G$ has the Chatterji-Ruane (CR) property.
\end{definition}

\begin{remark}\label{r:cr} Condition (i) is missing in \cite{CR}, but the authors of \cite{CR} informed me that it is used in the proof of the fact that the CR property implies RD and holds in all example of groups satisfying the CR property from \cite{CR}.
\end{remark}

\begin{theorem}\label{t2} The CR property implies the centroid property.
\end{theorem}

\proof We shall use $X$ as the centroid space of $G$. Let $\co$ be any map $G\times G\to X$ with the property $\co(g,k)\in C(x_0, g\cdot x_0) \cap C(x_0, k\cdot x_0) \cap C(g\cdot x_0, k\cdot x_0)$ where $x_0$ is a fixed base point (by (ii) such a map exists). We shall prove that $\co$ is a centroid map, that is conditions ($c_1$), ($c_2$) and ($c_3)$ from Definition \ref{d:0} hold.

($c_1$) Fix an element $k\in G$. Let $r>0$. Let $L(g)\le r$, $c(g)=\co(g,k)$. Then $c$ belongs to $C(x_0,g\cdot x_0)$. By (i) $C(x_0,g\cdot x_0)$ contains $x_0$. By (iv) the diameter of $C(x_0,g\cdot x_0)$ does not exceed $Q(r)$. By (iii) the number of possible $c(g)$ with $L(g)\le r$ does not exceed $R(Q(r))$. This gives ($c_1$).

($c_2$) is proved the same way as ($c_1$).

($c_3$) Fix an element $h\in G$. Consider all elements $g\in G$ with $L(g)\le r$. Then, by the definition,  $$\co(g,gh)\in  C(x_0, g\cdot x_0) \cap C(x_0, gh\cdot x_0) \cap C(g\cdot x_0, gh\cdot x_0.)$$
Therefore by the equivariance of the map $C$, we have
$$g\iv\cdot \co(g,gh)\in C(g\iv\cdot  x_0, x_0) \cap C(g\iv\cdot x_0, h\cdot x_0) \cap C(x_0, h\cdot x_0).$$ Hence $g\iv\cdot \co(g,gh)\in C(g\iv\cdot x_0,x_0)$. Since $\dist(g\iv\cdot x_0,x_0)\le r$, the set $C(g\iv \cdot x_0,x_0)$ is contained in the ball $B(g\iv \cdot x_0,Q(r))\le B(x_0,Q(r)+r)$ because $L(g)\le r$ (by conditions (i) and (iv)). Therefore the number of possible points of the form $g\iv\cdot \co(g,gh)$ does not exceed $R(Q(r)+r)$ (by (iii)).
\endproof

\begin{cor}[Chatterji-Ruane, \cite{CR}, Theorem 0.4] Groups acting properly with uniformly bounded stabilizers
and cellularly on a CAT(0) cube complex of finite dimension have the centroid property.\label{cor:ct}
\end{cor}

\begin{remark}
Corollary \ref{cor:ct} and Theorem \ref{th:mcg} have been recently generalized by Bowditch \cite{Bow14}. He proved that every coarsely median group satisfies property RD. \footnote{Note that the centroid ``cloud" map used in \cite{Bow14} is not necessarily equivariant as was noticed by Rudolf Zeidler. Moreover, even if it is equivariant, the argument \cite{Bow14} repeats the same mistake as \cite{DS} and \cite{BM}, see Footnote \ref{page}. Nevertheless, it is easy to deduce from \cite{Bow14} that these groups satisfy the centroid property.}For the precise precise definition of coarsely median groups see \cite[Page 170]{Bow14}. Informally, it means that that there exists a ``centroid" map $G^3\to G$ such that a) if we modify the triple $t\in G^3$ slightly, then the ``centroid" does not change much and b) every finite subset $A$ of $G$ is ``quasi-isomorphic" to a median algebra, with error constant depending only on the size of $A$. Both the groups acting ``nicely" on CAT(0) cube complexes and mapping class groups are coarsely median groups (see \cite{Bow14}).
\end{remark}

\begin{cor}[Osin, \cite{Os04}]\label{c:os} Every finitely generated group given by a (possibly infinite) set of relations satisfying the small cancelation property $C'(\lambda)$ with $\lambda<1/6$ has the CR property, and hence the centroid property and property RD.
\end{cor}

\proof Let $\Gamma$ be the (right) Cayley graph of $G$ corresponding to a finite generating set with the natural (left) action by $G$. For every pair of points $x,y\in \Gamma$ choose a geodesic path $[x,y]$. Then let $C(x,y)$ consist of all points on $[x,y]$ and all points on the loops $\gamma$ of $\Gamma$ labeled by the defining relations of $G$ such that $|\gamma\cap [x,y]|\ge \frac16|\gamma|$. The fact that this map satisfies conditions (i)-(iv) follows almost immediately from Strebel's description \cite{St} of geodesic triangles in the Cayley graphs of groups given by small cancelation presentations (see \cite{Os04}).
%
%
%
\endproof

\begin{remark}\label{rk} Property RD for groups given by presentations satisfying $C'(\frac{1}{10})$ was proved before by Arzhantseva and Dru\c tu \cite{AD}. \end{remark}

\subsection{The Ciobanu-Holt-Rees property}

Let $G$ be a finitely generated group with the word length function (with respect to some finite generating set).  We say that $(x,y)\in G\times G$ is a factorisation of $g\in G$ if $g=xy, L(x)+L(y)=L(g)$. Let $\mathcal{D}$ be a subset of the set of all decompositions of elements of $G$. Let $S(r)$ denote the set of all elements of length $r$ (the sphere of radius $r$ in the Cayley graph of $G$). For every $g\in S(r+r')$ let $\mathcal{F}_{g,r,r'}$ be the number of elements of the subset of $\mathcal{D}$ consisting of all decompositions $(x,y)\in \mathcal{D}$ of $G$ with $L(x)=r, L(y)=r'$. Let $\mathcal{F}_{\mathcal{D}, r,r'}$ be the supremum of all numbers $\F_{g,r,r'}$, $g\in S(r+r')$.

\begin{definition}[See \cite{CHR}]\label{d:chr} Suppose that

\begin{itemize}
\item[(D1)] $\F_{\D,r,r'}$ is bounded above by $P_1(\min(r, r'))$ for some polynomial $P_1(x)$.
\item[(D2)] For each $k \in G$, each $r, r'\in \R_+$, there is a subset $U(k, r, r')$ of $G\times G\times G$ as follows.
For each representation of $k$ as a product $g_1g_2$ with $g_1 \in S(r), g_2 \in S(r')$,
$U(k, r, r')$ contains a triple $(f_1,\hat g, f_2)$, for which $k = f_1\hat gf_2$, and $\hat k = h_1h_2$,
where $(f_1,h_1)\in \D$, $L(f_1)+L(h_1)=r,$ $f_1h_1=g_1$ and $(h_2, f_2)\in \D, L(h_2)+L(f_2)=r', h_2f_2=g_2$, $L(h_1),L(h_2)\le K\min(r,r')$ for some uniform constant $K$.
Furthermore, there are polynomials $P_2(x), P_3(x)$ such that
\begin{itemize}
\item[(a)] for all $k, r, r'$, $|U(k, r, r')| \le P_2(\min(r, r'))$,
\item[(b)] $|T (r, r')| \le P_3(\min(r, r'))$, where $T(r, r') = \{\hat k : \exists k, (f_1,\hat k, f_2) \in  U(k, r, r')\}.$
\end{itemize}
\end{itemize}
Then we say that $G$ satisfies the Ciobanu-Holt-Rees (CHR) property.
\end{definition}

\setlength{\unitlength}{0.70pt}
\begin{figure}
\begin{center}
{\large
\begin{picture}(500,260)(-250,-110)
\put(-250,0){\circle*{10}}
\put(0,-40){\circle*{10}}
\put(250,0){\circle*{10}}
\qbezier(-250,0)(0,50)(250,0) \put(0,30){$k$}
\qbezier(-250,0)(-125,0)(0,-40) \put(-140,-20){$g_1$}
\qbezier(250,0)(125,0)(0,-40) \put(125,-20){$g_2$}
\put(-100,-60){\circle*{10}}
\put(100,-60){\circle*{10}}
\qbezier(-250,0)(-175,-25)(-100,-60) \put(-160,-50){$f_1$}
\qbezier(250,0)(175,-25)(100,-60) \put(150,-50){$f_2$}
\qbezier(-100,-60)(0,-25)(100,-60) \put(0,-60){$\hat{k}$}
\put(-80,-45){$h_1$} \put(70,-45){$h_2$}
\end{picture}
}
\caption{\label{fig5} Condition D2 of the CHR property}
\end{center}
\end{figure}
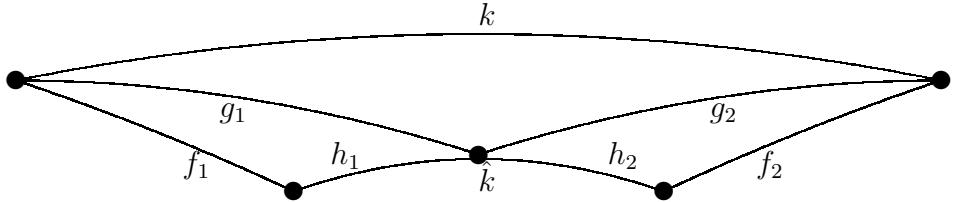
\begin{theorem}\label{t4} The CHR property implies the centroid property.
\end{theorem}

\proof Suppose that $G$ satisfies the CHR property with respect to a set of factorizations $\P$, polynomials $P_1, P_2, P_3$ and the word length function $L$. Let $\dist$ be the associated distance function: $\dist(g,h)=L(g\iv h)$. The group $G$ acts (on the left) on itself freely by isometries, so let the centroid space $X$ be $(G, \dist)$.

For every pair of elements $(g,k)$ of $G$ , $L(g)\le r$, let us pick one triple $(f_1,\hat g, f_2)$ from $U(k,r,L(g\iv k))$ as in (D2). Then let us define $\co(g,k)=f_1$. Let us prove that this function $\co$ satisfies parts  ($c_1$), ($c_2$) and ($c_3$) of the centroid property.

($c_1$)  Let $k\in G$, $r>0$. By part (a) of (D2)  there are at most $rP_2(r)$ triples in $U(k,r,L(g\iv k))$. Therefore there are at most $rP_2(r)$ choices for $\co(g,k)$.

($c_2$) Let $g\in G$. By (D1), there are at most $L(g)P_1(L(g))$ factorisations of $g$ from $\P$. Hence the there are at most $L(g)P_1(L(g))$ choices of $\co(g,k)$ ($k\in G$).

($c_3$) Let $h\in G$. Pick any $g\in G$. Let $k=gh$, so on Figure \ref{fig5} $g_1=g, g_2=h$. Let $f_1=\co(g,gh)$. Then $g\iv f_1=h_1$ where $(f_1,h_1)$ is a factorization of $g$ from $\P$. Let $h_1h_2=\hat k$. By (D2), the length of $h_2$ does not exceed $Kr$. Then the number of possible elements $h_2$ (with $h$ fixed) does not exceed $\F(D,Kr,L(h))\le P_1(Kr)$. The number of possible elements $\hat k=h_1h_2$ does not exceed $P_3(r)$ by part (b) of (D2). Therefore the number of possible elements $h_1=g\iv \co(g,gh)$ does not exceed $P_1(Kr)P_3(r).$
\endproof

Theorem \ref{t4} and \cite{CHR} imply the following

\begin{cor} Artin groups of large type have the centroid property.
\end{cor}

\subsection{The (*)-relative hyperbolicity}

\noindent \textit{Notation}: For a subset $Y$ in a metric space we
denote by $\overline{\nn}_\delta (Y)$ the closed tubular
$\delta$-neighborhood of $Y$, that is $\{ x\mid \dist (x,Y) \leq \delta
\}$.

\me

\begin{definition}\label{star}

Let $G$ be a group and let $H_1,...,H_m$ be subgroups in $G$. We say
that $G$ is (*)-{\textit{relatively hyperbolic with respect to}}
$H_1,...,H_m$ if there exists a finite generating set $S$ of $G$,
and two constants $\sigma$ and $\delta$ such that the following
property holds:

\begin{itemize}
\item[(*)] For every triple of points $A,B, C$ in the (right) Cayley graph of $G$ with resperct to $S$ pick geodesics connecting each pair of points: $\geod[A,B]$, $\geod[B,C]$ and  $\geod[A,C].$ Then there exists a coset $gH_i$ such that
$\overline{\nn}_\sigma(gH_i)$ intersects each of the geodesic sides of the
triangle, and the entrance (resp. exit) points $A_1, B_1, C_1$
(resp. $B_2, C_2, A_2$) of the sides $\geod[A_,B], \geod[B,C]$ and $\geod[C,A]$ in
$\overline{\nn}_\sigma(gH_i)$ satisfy $$\dist(A_1,A_2)<\delta,\,
\dist(B_1, B_2)<\delta,\, \dist(C_1, C_2)<\delta\, .$$
\end{itemize}
\end{definition}

\begin{figure}
\centering
\unitlength .7mm 
\linethickness{0.4pt}
\ifx\plotpoint\undefined\newsavebox{\plotpoint}\fi 
\begin{picture}(105.5,87)(0,0)

\qbezier(14.5,9.25)(42.38,38.38)(53.75,84)
\qbezier(53.75,84)(68.38,31.38)(101.5,10.25)
\qbezier(101.5,10.25)(59.5,24.25)(14.5,9.25)
\put(86.65,37.5){\line(0,1){1.26}}
\put(86.63,38.76){\line(0,1){1.258}}
\multiput(86.55,40.02)(-.0308,.3135){4}{\line(0,1){.3135}}
\multiput(86.43,41.27)(-.02872,.20804){6}{\line(0,1){.20804}}
\multiput(86.26,42.52)(-.03159,.17722){7}{\line(0,1){.17722}}
\multiput(86.03,43.76)(-.03369,.15386){8}{\line(0,1){.15386}}
\multiput(85.77,44.99)(-.03176,.12194){10}{\line(0,1){.12194}}
\multiput(85.45,46.21)(-.03319,.10964){11}{\line(0,1){.10964}}
\multiput(85.08,47.42)(-.0317,.0916){13}{\line(0,1){.0916}}
\multiput(84.67,48.61)(-.032746,.083838){14}{\line(0,1){.083838}}
\multiput(84.21,49.78)(-.033606,.076991){15}{\line(0,1){.076991}}
\multiput(83.71,50.94)(-.032291,.066719){17}{\line(0,1){.066719}}
\multiput(83.16,52.07)(-.032943,.061769){18}{\line(0,1){.061769}}
\multiput(82.57,53.18)(-.033479,.05725){19}{\line(0,1){.05725}}
\multiput(81.93,54.27)(-.032297,.050571){21}{\line(0,1){.050571}}
\multiput(81.25,55.33)(-.032697,.047027){22}{\line(0,1){.047027}}
\multiput(80.53,56.37)(-.033014,.043722){23}{\line(0,1){.043722}}
\multiput(79.77,57.37)(-.033256,.040629){24}{\line(0,1){.040629}}
\multiput(78.97,58.35)(-.033429,.037722){25}{\line(0,1){.037722}}
\multiput(78.14,59.29)(-.03354,.034984){26}{\line(0,1){.034984}}
\multiput(77.27,60.2)(-.034885,.033643){26}{\line(-1,0){.034885}}
\multiput(76.36,61.07)(-.037624,.03354){25}{\line(-1,0){.037624}}
\multiput(75.42,61.91)(-.040531,.033375){24}{\line(-1,0){.040531}}
\multiput(74.45,62.71)(-.043625,.033142){23}{\line(-1,0){.043625}}
\multiput(73.44,63.48)(-.046931,.032835){22}{\line(-1,0){.046931}}
\multiput(72.41,64.2)(-.050476,.032445){21}{\line(-1,0){.050476}}
\multiput(71.35,64.88)(-.057151,.033647){19}{\line(-1,0){.057151}}
\multiput(70.26,65.52)(-.061672,.033125){18}{\line(-1,0){.061672}}
\multiput(69.15,66.12)(-.066624,.032487){17}{\line(-1,0){.066624}}
\multiput(68.02,66.67)(-.072086,.031718){16}{\line(-1,0){.072086}}
\multiput(66.87,67.18)(-.083741,.032992){14}{\line(-1,0){.083741}}
\multiput(65.7,67.64)(-.091506,.031969){13}{\line(-1,0){.091506}}
\multiput(64.51,68.05)(-.10954,.03351){11}{\line(-1,0){.10954}}
\multiput(63.3,68.42)(-.12185,.03212){10}{\line(-1,0){.12185}}
\multiput(62.08,68.74)(-.13668,.03035){9}{\line(-1,0){.13668}}
\multiput(60.85,69.02)(-.17713,.03211){7}{\line(-1,0){.17713}}
\multiput(59.61,69.24)(-.20796,.02933){6}{\line(-1,0){.20796}}
\multiput(58.37,69.42)(-.3134,.0317){4}{\line(-1,0){.3134}}
\put(57.11,69.54){\line(-1,0){1.258}}
\put(55.85,69.62){\line(-1,0){2.52}}
\put(53.33,69.63){\line(-1,0){1.258}}
\put(52.08,69.56){\line(-1,0){1.254}}
\multiput(50.82,69.44)(-.24975,-.03373){5}{\line(-1,0){.24975}}
\multiput(49.57,69.27)(-.17731,-.03107){7}{\line(-1,0){.17731}}
\multiput(48.33,69.05)(-.15396,-.03324){8}{\line(-1,0){.15396}}
\multiput(47.1,68.79)(-.12203,-.0314){10}{\line(-1,0){.12203}}
\multiput(45.88,68.47)(-.10974,-.03287){11}{\line(-1,0){.10974}}
\multiput(44.67,68.11)(-.091692,-.03143){13}{\line(-1,0){.091692}}
\multiput(43.48,67.7)(-.083934,-.032499){14}{\line(-1,0){.083934}}
\multiput(42.31,67.25)(-.077089,-.033379){15}{\line(-1,0){.077089}}
\multiput(41.15,66.75)(-.066814,-.032095){17}{\line(-1,0){.066814}}
\multiput(40.01,66.2)(-.061866,-.032761){18}{\line(-1,0){.061866}}
\multiput(38.9,65.61)(-.057348,-.03331){19}{\line(-1,0){.057348}}
\multiput(37.81,64.98)(-.050666,-.032148){21}{\line(-1,0){.050666}}
\multiput(36.75,64.3)(-.047123,-.032558){22}{\line(-1,0){.047123}}
\multiput(35.71,63.59)(-.043819,-.032885){23}{\line(-1,0){.043819}}
\multiput(34.7,62.83)(-.040726,-.033136){24}{\line(-1,0){.040726}}
\multiput(33.72,62.04)(-.037821,-.033318){25}{\line(-1,0){.037821}}
\multiput(32.78,61.2)(-.035083,-.033437){26}{\line(-1,0){.035083}}
\multiput(31.87,60.33)(-.032496,-.033498){27}{\line(0,-1){.033498}}
\multiput(30.99,59.43)(-.033651,-.037525){25}{\line(0,-1){.037525}}
\multiput(30.15,58.49)(-.033494,-.040432){24}{\line(0,-1){.040432}}
\multiput(29.34,57.52)(-.03327,-.043527){23}{\line(0,-1){.043527}}
\multiput(28.58,56.52)(-.032972,-.046834){22}{\line(0,-1){.046834}}
\multiput(27.85,55.49)(-.032593,-.05038){21}{\line(0,-1){.05038}}
\multiput(27.17,54.43)(-.032124,-.0542){20}{\line(0,-1){.0542}}
\multiput(26.53,53.35)(-.033306,-.061574){18}{\line(0,-1){.061574}}
\multiput(25.93,52.24)(-.032683,-.066528){17}{\line(0,-1){.066528}}
\multiput(25.37,51.11)(-.031929,-.071992){16}{\line(0,-1){.071992}}
\multiput(24.86,49.96)(-.033238,-.083644){14}{\line(0,-1){.083644}}
\multiput(24.4,48.78)(-.032238,-.091412){13}{\line(0,-1){.091412}}
\multiput(23.98,47.6)(-.03102,-.10032){12}{\line(0,-1){.10032}}
\multiput(23.6,46.39)(-.03247,-.12175){10}{\line(0,-1){.12175}}
\multiput(23.28,45.18)(-.03075,-.13659){9}{\line(0,-1){.13659}}
\multiput(23,43.95)(-.03263,-.17703){7}{\line(0,-1){.17703}}
\multiput(22.77,42.71)(-.02994,-.20787){6}{\line(0,-1){.20787}}
\multiput(22.59,41.46)(-.0327,-.3133){4}{\line(0,-1){.3133}}
\put(22.46,40.21){\line(0,-1){1.257}}
\put(22.38,38.95){\line(0,-1){1.26}}
\put(22.35,37.69){\line(0,-1){1.26}}
\put(22.37,36.43){\line(0,-1){1.258}}
\put(22.43,35.17){\line(0,-1){1.255}}
\multiput(22.55,33.92)(.033,-.24985){5}{\line(0,-1){.24985}}
\multiput(22.72,32.67)(.03054,-.1774){7}{\line(0,-1){.1774}}
\multiput(22.93,31.42)(.03279,-.15406){8}{\line(0,-1){.15406}}
\multiput(23.19,30.19)(.03104,-.12213){10}{\line(0,-1){.12213}}
\multiput(23.5,28.97)(.03255,-.10983){11}{\line(0,-1){.10983}}
\multiput(23.86,27.76)(.031161,-.091784){13}{\line(0,-1){.091784}}
\multiput(24.26,26.57)(.032252,-.084029){14}{\line(0,-1){.084029}}
\multiput(24.72,25.39)(.033152,-.077187){15}{\line(0,-1){.077187}}
\multiput(25.21,24.24)(.031898,-.066908){17}{\line(0,-1){.066908}}
\multiput(25.76,23.1)(.032579,-.061962){18}{\line(0,-1){.061962}}
\multiput(26.34,21.98)(.033141,-.057446){19}{\line(0,-1){.057446}}
\multiput(26.97,20.89)(.033599,-.053298){20}{\line(0,-1){.053298}}
\multiput(27.64,19.83)(.032419,-.047219){22}{\line(0,-1){.047219}}
\multiput(28.36,18.79)(.032756,-.043916){23}{\line(0,-1){.043916}}
\multiput(29.11,17.78)(.033016,-.040823){24}{\line(0,-1){.040823}}
\multiput(29.9,16.8)(.033207,-.037918){25}{\line(0,-1){.037918}}
\multiput(30.73,15.85)(.033334,-.035181){26}{\line(0,-1){.035181}}
\multiput(31.6,14.93)(.033402,-.032594){27}{\line(1,0){.033402}}
\multiput(32.5,14.05)(.035987,-.032462){26}{\line(1,0){.035987}}
\multiput(33.44,13.21)(.040334,-.033613){24}{\line(1,0){.040334}}
\multiput(34.41,12.4)(.043429,-.033398){23}{\line(1,0){.043429}}
\multiput(35.4,11.64)(.046737,-.03311){22}{\line(1,0){.046737}}
\multiput(36.43,10.91)(.050284,-.032741){21}{\line(1,0){.050284}}
\multiput(37.49,10.22)(.054105,-.032283){20}{\line(1,0){.054105}}
\multiput(38.57,9.57)(.061476,-.033487){18}{\line(1,0){.061476}}
\multiput(39.68,8.97)(.066432,-.032879){17}{\line(1,0){.066432}}
\multiput(40.81,8.41)(.071898,-.032141){16}{\line(1,0){.071898}}
\multiput(41.96,7.9)(.083546,-.033484){14}{\line(1,0){.083546}}
\multiput(43.13,7.43)(.091316,-.032507){13}{\line(1,0){.091316}}
\multiput(44.31,7.01)(.10023,-.03131){12}{\line(1,0){.10023}}
\multiput(45.52,6.63)(.12166,-.03283){10}{\line(1,0){.12166}}
\multiput(46.73,6.3)(.1365,-.03116){9}{\line(1,0){.1365}}
\multiput(47.96,6.02)(.17693,-.03315){7}{\line(1,0){.17693}}
\multiput(49.2,5.79)(.20778,-.03056){6}{\line(1,0){.20778}}
\multiput(50.45,5.61)(.3132,-.0336){4}{\line(1,0){.3132}}
\put(51.7,5.47){\line(1,0){1.257}}
\put(52.96,5.39){\line(1,0){1.26}}
\put(54.22,5.35){\line(1,0){1.26}}
\put(55.48,5.36){\line(1,0){1.259}}
\put(56.73,5.43){\line(1,0){1.255}}
\multiput(57.99,5.54)(.24994,.03226){5}{\line(1,0){.24994}}
\multiput(59.24,5.7)(.17749,.03002){7}{\line(1,0){.17749}}
\multiput(60.48,5.91)(.15416,.03234){8}{\line(1,0){.15416}}
\multiput(61.72,6.17)(.12222,.03068){10}{\line(1,0){.12222}}
\multiput(62.94,6.48)(.10993,.03222){11}{\line(1,0){.10993}}
\multiput(64.15,6.83)(.09953,.03346){12}{\line(1,0){.09953}}
\multiput(65.34,7.23)(.084123,.032005){14}{\line(1,0){.084123}}
\multiput(66.52,7.68)(.077284,.032925){15}{\line(1,0){.077284}}
\multiput(67.68,8.17)(.071189,.033683){16}{\line(1,0){.071189}}
\multiput(68.82,8.71)(.062057,.032397){18}{\line(1,0){.062057}}
\multiput(69.93,9.3)(.057543,.032972){19}{\line(1,0){.057543}}
\multiput(71.03,9.92)(.053397,.033442){20}{\line(1,0){.053397}}
\multiput(72.1,10.59)(.047314,.03228){22}{\line(1,0){.047314}}
\multiput(73.14,11.3)(.044012,.032627){23}{\line(1,0){.044012}}
\multiput(74.15,12.05)(.04092,.032896){24}{\line(1,0){.04092}}
\multiput(75.13,12.84)(.038016,.033095){25}{\line(1,0){.038016}}
\multiput(76.08,13.67)(.035279,.03323){26}{\line(1,0){.035279}}
\multiput(77,14.53)(.032692,.033306){27}{\line(0,1){.033306}}
\multiput(77.88,15.43)(.032568,.035891){26}{\line(0,1){.035891}}
\multiput(78.73,16.37)(.033731,.040235){24}{\line(0,1){.040235}}
\multiput(79.54,17.33)(.033526,.043331){23}{\line(0,1){.043331}}
\multiput(80.31,18.33)(.033247,.046639){22}{\line(0,1){.046639}}
\multiput(81.04,19.35)(.032889,.050188){21}{\line(0,1){.050188}}
\multiput(81.73,20.41)(.032442,.05401){20}{\line(0,1){.05401}}
\multiput(82.38,21.49)(.033667,.061377){18}{\line(0,1){.061377}}
\multiput(82.99,22.59)(.033074,.066335){17}{\line(0,1){.066335}}
\multiput(83.55,23.72)(.032352,.071803){16}{\line(0,1){.071803}}
\multiput(84.07,24.87)(.03373,.083447){14}{\line(0,1){.083447}}
\multiput(84.54,26.04)(.032775,.09122){13}{\line(0,1){.09122}}
\multiput(84.96,27.22)(.03161,.10014){12}{\line(0,1){.10014}}
\multiput(85.34,28.43)(.03319,.12156){10}{\line(0,1){.12156}}
\multiput(85.67,29.64)(.03156,.13641){9}{\line(0,1){.13641}}
\multiput(85.96,30.87)(.03367,.17684){7}{\line(0,1){.17684}}
\multiput(86.19,32.11)(.03117,.20769){6}{\line(0,1){.20769}}
\multiput(86.38,33.35)(.0276,.2505){5}{\line(0,1){.2505}}
\put(86.52,34.61){\line(0,1){1.257}}
\put(86.61,35.86){\line(0,1){1.638}}
\put(12,7.75){\makebox(0,0)[cc]{$A$}}
\put(54.25,87){\makebox(0,0)[cc]{$B$}}
\put(105.5,8.25){\makebox(0,0)[cc]{$C$}}
\put(23,22.25){\makebox(0,0)[cc]{$A_1$}}
\put(61,73.5){\makebox(0,0)[cc]{$B_1$}}
\put(78.5,10){\makebox(0,0)[cc]{$C_1$}}
\put(46.25,73.25){\makebox(0,0)[cc]{$B_2$}}
\put(89.5,25.75){\makebox(0,0)[cc]{$C_2$}}
\put(31.25,10){\makebox(0,0)[cc]{$A_2$}}
\put(49.5,69.25){\circle*{1.41}} \put(58.5,69.5){\circle*{1.8}}
\put(84.25,24.75){\circle*{1.12}} \put(77.75,15.75){\circle*{1.41}}
\put(32.75,14){\circle*{1.5}} \put(26,22.5){\circle*{1}}
\put(88.5,54){\makebox(0,0)[cc]{$\overline{N}_\sigma(gH_i)$}}
\end{picture}
\caption{\label{fig2} (*)-relative hyperbolicity}
\end{figure}
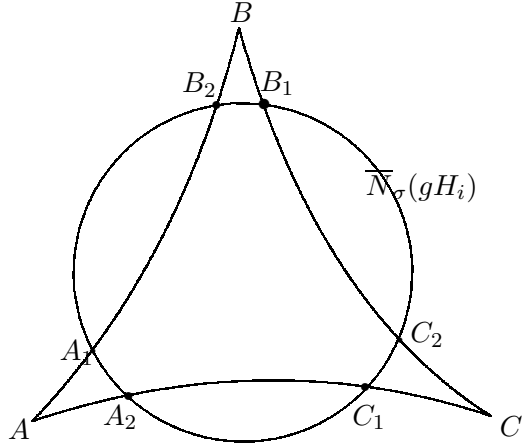

By  \cite[Proposition 2.9]{DS} if $G$ is (strongly) relatively hyperbolic with respect to subgroups $H_1,\ldots ,H_m$, then it is (*)-relatively hyperbolic with respect to these subgroups.

\begin{theorem} \label{t:*} If $G$ is (*)-relatively hyperbolic with respect to a subgroup $H$, and $H$ satisfies the centroid property, then $G$ satisfies the centroid property.
\end{theorem}

\proof Let $X$ be a centroid space for $H$ and $\co_H$ be its centroid map. We assume that the action of $G$ on $X$ is free. The general case of almost free actions is similar.

We construct a centroid space for $G$ as a representation induced by the representation of $H$ (in $\mathrm{Iso}(X))$ in the following fairly standard manner. Let $Y=X\times G$.
Define an action of $G$ on $Y$ by $g\cdot (x,k)=(x,g\iv k)$.
Now define an equivalence relation $\sim$ on $Y$ by
$(x,k)\sim (h\cdot x, kh)$ for every $h\in H$. Let $Z=Y/\si$. Note that $\sim$ respects the action of $G$ on $Y$: for every $g,k\in G, h\in H, x\in X$, $g\cdot (x,k)=(x,g\iv k)\sim (h\cdot x, g\iv kh)$.  Therefore $G$ acts in a natural way on $Z$.

The set $X$ maps into $Z$ by $x\mapsto (x,1)/\si$. It is easy to see that this map is injective,
and is equivariant with respect to the action of $H$. The action of $G$ on $Z$ is free. Indeed, if $g\cdot (x,k)=(x,g\iv k)\sim (x,k)$ for some $g,k,x$, then there exists $h\in H$ such that $h\cdot x=x, g\iv kh=k$. Since the action of $H$ on $X$ is free, we deduce that $h=1$, hence $g=1$.

Let $U$ be a set of representatives of all left cosets $gH$ of the subgroup $H$ in $G$.
We assume that for every coset $tH$, its representative in $U$ has the smallest possible length in $G$. In particular, $1$ is the representative of the coset $H$.
It is easy to check that every $\sim$-class contains unique element of the form $(x,t)$, $t\in U$. Thus we can identify $Z$ with the set $Z'=\{(x,t)\mid x\in X, t\in U\}$. The corresponding action of $G$ on $Z'$ is defined as follows: $g(x,t)=(h\iv \cdot x, t')$ where $g\iv t=t'h, t'\in U, h\in H$.
Let us define a metric $\dist$ on $Z$ by
$\dist((x,1),(y,t))=\dist(x,y)+L(t)$ and $\dist((x,t), (y,t'))=\dist((x,1),(h\cdot y,t'')$ where $h\in H, t''\in U$ and $t\iv t' =t''h$. It is  straightforward to verify that the restriction of $\dist$ on $X=X\times \{1\}$ coincides with the metric of $X$, and that $G$ acts on the metric space $(Z,\dist)$ by isometries.

Now let us define the centroid map. Let $g,k\in G$. Consider a geodesic triangle $(1,g,k)$ in the Cayley graph of $G$. Denote its geodesic sides
by $\geod[1,g]$, $\geod[g,h]$, $\geod[1,k]$ (pick one geodesic for each pair of vertices of the triangle).
Then there exists a coset $tH$, $t\in U$, whose $\delta$-neighborhood satisfies condition (*). This means that there exist three elements $u,h_g, h_k\in H$ such that the following conditions hold.

\begin{itemize}
\item[(1)] $tuh_g$ is at distance at most $\sigma+2\delta$ from the geodesic sides $\geod[1,g]$ and $\geod[g,k]$.

\item[(2)] $tuh_k$ is at distance at most $\sigma+2\delta$ from the geodesic sides $\geod[1,k]$ and $\geod[g,k]$.

\item[(3)] $tu$ is at distance at most $\sigma+2\delta$ from the geodesic sides $\geod[1,g]$ and $\geod[1,k]$.
\end{itemize}

Then let us define the centroid function $\co_G\colon G\times  G\to Z$ by $\co_G(g,k)=\co(uh_g, uh_k) \times t$, that is, we add $t$ as the second coordinate to each vertex of $\co(uh_g, uh_k)$ to obtain an element from $Z'$.

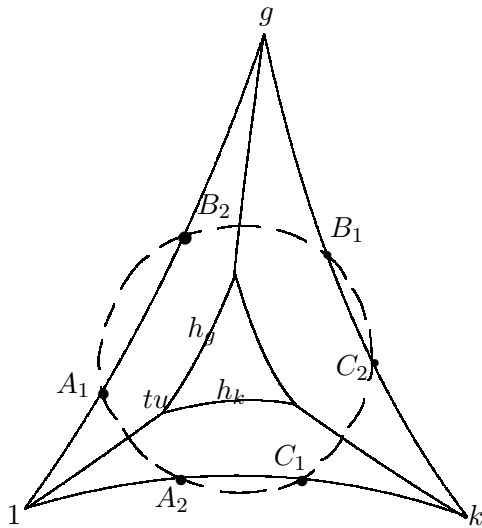
\begin{figure}[!ht]
\centering
\unitlength .7mm 
\linethickness{0.4pt}
\ifx\plotpoint\undefined\newsavebox{\plotpoint}\fi 
\begin{picture}(91.25,93.75)(0,20)
\qbezier(7.25,23.5)(35,61.38)(52.75,113.75)
\qbezier(52.75,113.75)(65,56.13)(91.25,22)
\qbezier(91.25,22)(50.25,36.75)(7.25,23.5)
\qbezier(7.5,23.5)(20,31.88)(33.5,41.75)
\qbezier(33.5,41.75)(41.13,52.13)(47.25,68)
\qbezier(47.25,68)(52.5,49.88)(58.75,43.25)
\qbezier(58.75,43.25)(47.75,45.5)(33.75,41.75)
\qbezier(58.25,43.5)(74.38,31.88)(91,21.75)
\qbezier(47,67.5)(48.25,80.88)(52.5,112.75)
\put(46.5,45.75){\makebox(0,0)[cc]{$h_k$}}
\put(41,57.25){\makebox(0,0)[cc]{$h_g$}}
\qbezier(25.5,65.25)(28.38,70.13)(30.75,71.5)
\qbezier(33,73.25)(36.13,75.25)(40.75,76.25)
\qbezier(44,76.75)(48.5,77.63)(52,77)
\qbezier(55.25,76.5)(58.5,75.63)(60.75,74.25)
\qbezier(63.5,72.5)(67.13,70.13)(68.25,67.25)
\qbezier(69.5,65.25)(71.25,63.13)(72,59.5)
\qbezier(72.75,57.5)(73.88,48.25)(72.5,48)
\qbezier(71.75,44.25)(71.25,41.88)(68.75,39)
\qbezier(67.75,37.25)(65.38,33)(63.5,32.75)
\qbezier(62.25,30.75)(59.63,28.5)(55.5,27.25)
\qbezier(53,27)(48.75,26.13)(44.5,26.75)
\qbezier(41,27.75)(35.88,29.25)(32.25,31.75)
\qbezier(30.5,33)(28.25,34.75)(26,38.5)
\qbezier(24.75,40.5)(22.25,43)(21.75,46.5)
\qbezier(21.5,50)(20.75,51.25)(22,55.5)
\qbezier(22.5,58.25)(23.63,60.75)(24.25,62.25)
\put(16.5,47){\makebox(0,0)[cc]{$A_1$}}
\put(43.25,81){\makebox(0,0)[cc]{$B_2$}}
\put(53.25,117){\makebox(0,0)[cc]{$g$}}
\put(68.5,76.5){\makebox(0,0)[cc]{$B_1$}}
\put(69.5,50.25){\makebox(0,0)[cc]{$C_2$}}
\put(22.25,45.25){\circle*{2.06}} \put(37.75,75){\circle*{2.5}}
\put(64.75,71.75){\circle*{1.58}} \put(73.75,51.25){\circle*{1.58}}
\put(60,28.75){\circle*{2.24}} \put(37,29){\circle*{2.12}}
\put(35.25,25.25){\makebox(0,0)[cc]{$A_2$}}
\put(5.25,22.25){\makebox(0,0)[cc]{$1$}}
\put(57.75,33.25){\makebox(0,0)[cc]{$C_1$}}
\put(93,22.25){\makebox(0,0)[cc]{$k$}}
\put(32.25,44.25){\makebox(0,0)[cc]{$tu$}}
\end{picture}
\caption{\label{fig3} Triangle from $H$ inscribed in a triangle from $G$}
\end{figure}

Property ($c_1$) is proved as follows. Fix an element $k\in G$ and $r>0$,  and let us vary $g\in G, L(g)\le r$. Then since $tu$ is at distance at most $\sigma+2\delta$ from $\geod[1,g]$ and from $\geod[1,k]$, it is at distance at most $\sigma+2\delta$ from a point $A_2$ on $\geod[1,k]$ (see Figure \ref{fig3} which is at distance at most $r+4(\sigma+\delta)$ from 1. The number of such points $A_2$ is at most $r+2(\sigma+2\delta)$. The number of elements inside a ball of radius $2(\sigma+2\delta)$ is a universal constant. Therefore the number of possibilities for $tu$ is at most $K_1r+K_2$ for some constants $K_1, K_2$. Since $t$ is the only element from $U$ in the coset $tuH$, the number of possibilities for $t$ and the number of possibilities for $u$ do not exceed $K_1r+K_2$.  Since the point $C_1$ is determined uniquely by the coset $tH$, we have at most $K_1r+K_2$ possibilities for the point $C_1$. The distance from $tuh_k$ to $C_2$ is at most $\sigma+2\delta$. Hence the number of possibilities for $tuh_k$ is at most $K_3r+K_4$ for some constants $K_3, K_4$. Thus the number of possible choices for $h_k$ does not exceed $P_2(r)$ where $P_2$ is a (universal) polynomial of degree 2. The length of $h_g$ does not exceed $r+2(\sigma+2\delta)$. Therefore by the centroid property (for $H$), the number of possibilities for $\co_H(uh_g,uh_k)$ does not exceed $P(r)$ where $P$ is a (universal) polynomial. Therefore the number of possibilities for $\co_G(g,k)$ does not exceed the number of possibilities for $t$ times $P(r)$, this is can be bounded from above by $(K_1r+K_2)P(r)$.

Properties ($c_2$) and ($c_3$) are proved in a similar way and we are leaving the checking as an exercise for the reader.
\endproof

\begin{remark} Theorem \ref{t:*} can be easily generalized to groups that are (*)-relatively hyperbolic with respect to several subgroups $H_1,\ldots, H_m$ satisfying the centroid property.
\end{remark}

The proof of the following statement is very similar to (but easier than) the proof of Theorem \ref{t:*}, and is left to the reader.


\begin{theorem}\label{t:ab} (a) Suppose that $G$ is (*)-relatively hyperbolic with respect to subgroups $H_1,\ldots ,H_m$. Then $G$ has relative centroid property with respect to subgroups $H_i$, $i=1,\ldots,m$ (see Definition \ref{de:rd}). Moreover if each $H_i$ has $RD$, then each $T_i=G\cdot H_i^3$ has RD.

(b) Suppose that $G$  has the  relative centroid property with respect to subgroups $H_1,\ldots,H_m$ each of which has the centroid property. Then $G$ has the centroid property.
\end{theorem}

%

\subsection{Lafforgue's properties $(H_\delta)$ and $(K_\delta)$}

Conditions $(H_\delta)$ and $(K_\delta)$ were used by Lafforgue in \cite{L}. Similar conditions (for $\delta=0$) were used by Ramagge, Robertson and Steger in \cite{RRS}.

Let $X$ be a metric space, $\delta\ge 0$. We say that a sequence of points $x_1x_2\ldots x_n$ in $X$ is a $\delta$-path if $$\dist(x_1,x_2)+\ldots+\dist(x_{n-1},x_n)\le \dist(x_1,x_n)+\delta.$$

\begin{definition}[Property $(H_\delta)$] Let $\delta\geq 0$. A discrete metric space $(X,d)$ satisfies \emph{property $(H_\delta)$} if there exists a polynomial $P_\delta$ such that for any $r\in{\R}_+$, $x,y\in X$ the set
$$\{t\in X\hbox{ such that } xty \hbox{ is a } \delta\hbox{-path, } \dist(x,t)\le r\}$$
contains at most $P_\delta(r)$ elements.

\end{definition}


We say that a triple $(x,y,z)\in X^3$ is $\delta$-\emph{retractibe} if there exists a point $t\in X$ such that $xty$, $ytz$ and $ztx$ are $\delta$-paths. (So if $X$ is hyperbolic, then every triple is $\delta$-retractible where $\delta$ depends only on the hyperbolicity constant.)

\begin{definition}[Property $(K_\delta)$]
Let $(X,d)$ be a metric space and $\Gamma$ be a discrete group acting almost freely by isometries on $X$, $\delta\ge 0$. We say that the pair $(X,\Gamma)$ satisfies {\it property $(K_\delta)$} if there exists $k\in \N$ and $\Gamma$-invariant subsets ${\mathcal T}_1,\dots,{\mathcal T}_k$ of $X^3$ such that:
\begin{itemize}\item[$(K_\delta a)$] There exists $C\in \R_+$ and a map $\lambda\colon X^3\to \cup_{i=1}^m T_i$ such that if $$\lambda(x,y,z)=(\alpha,\beta,\gamma),$$ then
\begin{equation}\label{101}\max\{d(\alpha,\beta),d(\beta,\gamma),d(\gamma,\alpha)\}\leq C \min\{d(x,y),d(y,z),d(z,x)\}+\delta\end{equation}
and $x\alpha\beta y$, $y\beta\gamma z$, $z\gamma\alpha x$ are $\delta$-paths.
\medskip
\item[$(K_\delta b)$] For any $i\in\{1,\dots,k\}$ and $\alpha,\beta,\gamma,\gamma'\in X$, if $(\alpha,\beta,\gamma)\in{\mathcal T}_i$ and $(\alpha,\beta,\gamma')\in{\mathcal T}_i$ then the triples $(\alpha,\gamma,\gamma')$ and $(\beta,\gamma,\gamma')$ are $\delta$-retractable.\end{itemize}
\end{definition}

\begin{figure}
\begin{center}
\unitlength .7mm 
\linethickness{0.4pt}
\ifx\plotpoint\undefined\newsavebox{\plotpoint}\fi 
\begin{picture}(94.5,97.25)(0,20)
\qbezier(7.25,23.5)(35,61.38)(52.75,113.75)
\qbezier(52.75,113.75)(65,56.13)(91.25,22)
\qbezier(91.25,22)(50.25,36.75)(7.25,23.5)
\qbezier(7.5,23.5)(20,31.88)(33.5,41.75)
\qbezier(33.5,41.75)(41.13,52.13)(47.25,68)
\qbezier(47.25,68)(52.5,49.88)(58.75,43.25)
\qbezier(58.75,43.25)(47.75,45.5)(33.75,41.75)
\qbezier(58.25,43.5)(74.38,31.88)(91,21.75)
\qbezier(47,67.5)(48.25,80.88)(52.5,112.75)
\put(5.5,23){\makebox(0,0)[cc]{$x$}}
\put(53.75,117.25){\makebox(0,0)[cc]{$y$}}
\put(94.5,22){\makebox(0,0)[cc]{$z$}}
\put(33.5,38.75){\makebox(0,0)[cc]{$\alpha$}}
\put(49.75,69){\makebox(0,0)[cc]{$\beta$}}
\put(61,44.25){\makebox(0,0)[cc]{$\gamma$}}
\end{picture}
\caption{\label{fig7} Condition $(K_\delta a)$.}
\end{center}

\end{figure}
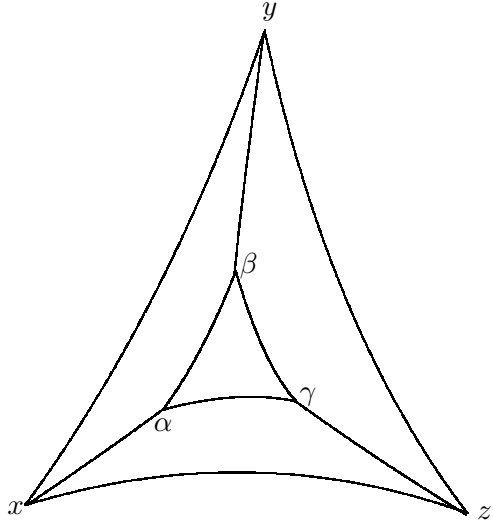

It was proved in \cite{L} that if a group acts almost freely on a metric space $X$, and the pair $(X,G)$ satisfies the conditions $(H_\delta)$ and $(K_\delta)$ for some $\delta$, then $G$ has property RD.

\begin{prob} \label{p09} Do conditions $(H_\delta)$ and $(K_\delta)$ imply the centroid property?
\end{prob}

Although we do not know the answer to this problem, the following statement holds.

\begin{theorem}\label{tl} If a group $G$ satisfies $(H_\delta)$ and $(K_\delta a)$ for some $\delta$ and some $G$-invariant sets $T_1,...,T_m\subseteq X^3$, then $G$ has the relative centroid property with respect to $T_1,\ldots, T_m$.
\end{theorem}
\proof Let $\lambda$ be the map from $(K_\delta a)$. Pick a point $x_0\in X$. For every $g,k\in G$ let $\rc(g,k)=\lambda(x_0,g\cdot x, k\cdot x)$,We will show that $\rc$ satisfies ($rc_1$)-($rc_4$).

($rc_1$) Pick $k\in G$ and let us vary $g, L(g)\le r$. Let $(\alpha,\beta,\gamma)=\rc(g,k)=\lambda(x_0,g\cdot x_0, g\cdot x_0)$. Then by $(K_\delta a)$ $x_0\alpha\beta(g\cdot x_0)$ is a $\delta$-path. Hence $\dist(x_0,\alpha)\le r+\delta$ because $L(g)\le r$. Since $(k\cdot x_0)\alpha x_0$  is a $\delta$-path, by $(H_\delta)$ the number of possible points $\alpha$ does not exceed $P_\delta(r+\delta)$. Since $\dist(\alpha,\gamma)\le Cr+\delta$ by (\ref{101}), and $\alpha\gamma(k\cdot x_0)$ is a $\delta$-path, by $(H_\delta)$ the number of possible points $\gamma$ for a given $\alpha$ does not exceed $P_\delta(Cr+\delta)$. Thus the possible number of pairs $(\alpha,\beta)$ does not exceed $P_\delta(r+\delta)P_\delta(Cr+\delta)$.

($rc_2$) is proved in the same way as ($rc_1$).

($rc_3$) Pick $h\in G$ and let us vary $g\in G$, $L(g)\le r$. Let $(\alpha,\beta,\gamma)=\rc(g,gh)=\lambda(x_0,g\cdot x_0, (gh)\cdot x_0)$. As in the proof of ($rc_1$), $\dist(g\cdot x_0, \beta)\le r+\delta$, hence $\dist(x_0,g\iv\cdot \beta)\le r$. Since $x_0(g\iv\cdot \beta)(h\cdot x_0)$ is a $\delta$-path, the number of possible points $g\iv \beta$ does not exceed $P_\delta(r+\delta)$. Then, as before, given $\beta$, the number of possible points $g\iv \cdot \gamma$ does not exceed $P_\delta(Cr+\delta)$. Therefore the number of possible pairs $g\iv \cdot (\beta,\gamma)$ does not exceed $P_\delta(r+\delta)P_\delta(Cr+\delta)$.

($rc_4$) immediately follows from (\ref{101}).
\endproof

\begin{remark} By \cite[Proposition 2.3]{L}, every set of triples $T\subseteq X^3$ satisfying $H_\delta$ and $(K_\delta b)$ has property RD. Thus properties $(H_\delta)$ and $(K_\delta)$ imply RD by Theorem \ref{t10} (of course that result is proved in \cite{L} too\footnote{Note, though, that \cite[Lemma 3.6]{L} which is supposed to define the relative centroid map does not have a proof (and the map is not constructed in \cite{L}). A proof is substituted by a not very precise reference to \cite{RRS}.}). It would be interesting (in view of Problem \ref{p09}) to find out if $(H_\delta)$, $(K_\delta b)$ imply a (correctly formulated) centroid property for each $T_i$.
\end{remark}

\subsection{Graph products of groups}

\subsubsection{Direct products}
As far as I know the only proof of property RD that does not involve a centroid-like property is the proof by Jolissaint \cite{J} that a direct product of two groups $A, B$ (and more general ``polynomially growing" extensions of $A$ by $B$) have RD if and only if both $A$ and $B$ have RD.
For completeness we present the proof for direct products here. It is based on the proof from \cite{J} and the clarifications sent to us by Paul Jolissaint.

\begin{theorem}[Jolissaint \cite{J}]\label{t:j} Let $A, B$ be two countable groups with length functions $L_A, L_B$, let $A\times B$ be the direct product with length function $L(a,b)=L(a)+L(b)$. Then $A\times B$ has property RD with respect to $L$ if and only if $A$ has property RD with respect to $L_A$, and $B$ has property RD with respect to $L_B$.
\end{theorem}

\proof The ``only if" part is obvious because the restriction of $L$ on $A$ (resp. $B$) coincides with $L_A$ (resp. $L_B$). Suppose that both $A$ and $B$ have property RD, and $P_A, P_B$ are the corresponding polynomials. Let $\phi, \psi$ be two functions $A\times B\to \R_+$ with finite supports and $\prop(\phi)=r$. In what follows $a, \alpha$ denote elements from $A$, $b, \beta$ denote elements from $B$. We need to estimate
\begin{equation}\label{01}||\phi*\psi||^2=\sum\limits_{a,b}\left(\sum\limits_{\alpha,\beta} \phi(\alpha,\beta)\psi(\alpha\iv a, \beta\iv b)\right)^2. \end{equation}

Let us denote $\phi(\alpha,\beta)$ by $\phi_\alpha(\beta)$ and $\psi(\alpha, \beta)$ by $\psi_{\alpha}(\beta)$. That is, for each $\alpha\in A$ we introduce two functions $\phi_\alpha, \psi_\alpha\colon B\to \R_+$. Then (\ref{01}) can be rewritten as

\begin{equation}\label{02}
\sum\limits_a\sum\limits_b\left(\sum\limits_{\alpha}\sum\limits_{\beta}\phi_\alpha(\beta)\psi_{\alpha\iv a}(\beta\iv b)\right)^2.
\end{equation}

For every $a\in A, b\in B$ let us denote $\sum\limits_{\alpha}\sum\limits_{\beta}\phi_\alpha(\beta)\psi_{\alpha\iv a}(\beta\iv b)$ by $f_a(b)$. That is, we introduce a new function $f_a\colon B\to \R_+$ for each $a\in A$. Then (\ref{02}) can be rewritten as
\begin{equation}\label{03}
\sum\limits_a\sum\limits_b f_a(b)^2=\sum\limits_a||f_a||^2.
\end{equation}
Note that $f_a=\sum\limits_\alpha\phi_\alpha*\psi_{\alpha\iv a}$. Therefore (\ref{03}) can be rewritten as
$$
\sum\limits_a\left|\left|\sum\limits_\alpha\phi_\alpha*\psi_{\alpha\iv a}\right|\right|^2
$$
which does not exceed

\begin{equation}\label{04}
\sum\limits_a\left(\sum\limits_\alpha\left|\left|\phi_\alpha*\psi_{\alpha\iv a}\right|\right|\right)^2
\end{equation}
by the triangle inequality.

Note that for each $\alpha$, $\prop(\phi_\alpha)\le r$. Therefore by property RD for $B$ (\ref{04}) can be estimated from above by
\begin{equation}\label{05}
P_B(r)\sum\limits_a\left(\sum\limits_\alpha||\phi_\alpha||\cdot ||\psi_{\alpha\iv a}||\right)^2.
\end{equation}

Now for every $a\in A$ let $\Phi(a)=||\phi_a||, \Psi(a)=||\psi_a||$. Thus we introduced two functions $\Phi, \Psi\colon A\to \R_+$. It is easy to check that
$$||\Phi||=||\phi||, ||\Psi||=||\psi||.$$

Then (\ref{05}) can be rewritten as

$$P_B(r)||\Phi*\Psi||^2.$$

Since the propagation of $\Phi$ does not exceed $r$, by property RD for $A$ we can estimate this from above by

$$P_A(r)P_B(r)||\Phi||^2||\Psi||^2$$
which is equal to $P_A(r)P_B(r)||\phi||^2||\psi||^2.$
\endproof

\begin{remark} It would be interesting to find a property which would generalize the relative centroid property and hold for direct products. Perhaps a multi-dimensional version of the relative centroid property should be defined for this purpose. That could help dealing with uniform lattices in semi-simple Lie groups of higher ranks.
\end{remark}

\begin{remark} \label{r:d} It is easy to check that if a group $A$ (resp. $B$) has the relative centroid property with space of centroids $X$, sets of triples $T_1,
\ldots, T_m$, and a relative centroid map $\rc_A$ (resp. $Y, T_1',\ldots,T_n', \rc_B$), then $A\times B$ has the relative centroid property with respect to $X\times Y$, $T_i\times T_j, 1\le i\le m, 1\le i\le n$ and $\rc_A\times \rc_B\mapsto ((a,b), (\alpha,\beta))\to (\rc_A(a,b), \rc_B(\alpha,\beta))$ (where we identify $(X\times Y)^3$ with $X^3\times Y^3$ in the natural way). In particular if both $A$ and $B$ have the centroid property, then $A\times B$ has the centroid property.
\end{remark}

\subsubsection[Graph products]{Graph products\footnote{This subsection is written jointly with Mitchel Kleban as a part of his Summer 2014 REU project.}}

Graph products of groups generalize both direct products and free products. Let $\Gamma=(V,E)$ be a finite unoriented graph with vertex set $V$ and edge set $E$ without loops and multiple edges. For each $v\in V$ let $G_v$ be a group. The graph product $G=\prod_\Gamma G_v$ is the quotient of the free product $\bigast\limits_{v\in V} G_v$ by the normal subgroup generated by all commutators $[g_1,g_2]$ where $g_1\in G_{v_1}, g_2\in G_{v_2}, (v_1,v_2)\in E$. Every element $g\in G$ is a product of \emph{syllables}
\begin{equation}\label{11}
g=g_{v_1}\ldots g_{v_m},
\end{equation} where $g_{v_i}\in G_{v_i}, v_i\ne v_{i+1}$ for every $i$. The minimal such $m$ is called the \emph{syllable length} of $g$, denoted by $\lambda(g)$. It is known \cite[Theorem 3.9]{rr} that any two representations (\ref{11}) of $g$ of minimal length (such representations are called \emph{reduced}) differ only by the order of the syllables. Thus if we have a length function $L_v$ on each group $G_v$, we can define a length function $L$ on $G$ by $$L(g)=\sum\limits_{i=1}^m L_{v_i}(g_i)+\lambda(g)$$ for every minimal representation (\ref{11}). The length function induced by $L$ on each $G_v$ is equivalent to $L_v$ because $\lambda(g)=1$ for each $g\in G_v$. A representation $g=u_1u_2\ldots u_n$ will be called a \emph{factorization} if $\lambda(g)=\lambda(u_1)+\ldots+\lambda(u_n)$.

For every full subgraph $\Gamma'=(V',E')\subseteq \Gamma$ let $G_{\Gamma'}$ be the corresponding graph product $\prod_{\Gamma'} G_v$. It is clear \cite{rr} that $G_{\Gamma'}$ is a subgroup of $G$ (and the natural map from $G_{\Gamma'}$ to $G$ is injective).

Let $\cc$ be the set of all cliques of $\Gamma$ (that is, subsets of $V$ where each pair of vertices is connected by an edge). For every $C\in\cc$ we call $G_{C}$ a \emph{clique subgroup} of $G$. Clearly every clique subgroup $G_C$ is the direct product of the vertex groups $G_v, v\in C$. Moreover since the syllable length of every element of a clique subgroup $G_C$ does not exceed $|C|$, the restriction of the length function $L$ of $G$ to $G_C$ is equivalent to the natural length function of $G_C:$ $$L_C(g_{v_1}\ldots g_{v_m})=L_{v_1}(g_{v_1})+\ldots+L_{v_m}(g_{v_m})$$
where $g_{v_i}\in G_{v_i}, v_i\in C$, $i=1,\ldots, m$.

The main result of \cite{CHR1} is that the graph product $G$ has property RD with respect to the length function $L$ if and only if each $G_v, v\in V,$ has property RD with respect to the length function $L_v$. This fact follows from the next theorem, Theorem \ref{t:j} and Theorem \ref{t10}. Note though that the proof of the next theorem is based on the intermediate results from \cite{CHR1}.

\begin{theorem} \label{t:c} The graph product $G$ has the relative centroid property with respect to the clique subgroups $G_C$, $C\in \cc$ (see Definition \ref{de:rd}).
\end{theorem}

\proof Let $g,h\in G$, $k=gh$ with $\lambda(k)=\lambda(g)+\lambda(h)-q$, $q\ge 0$. Then by \cite[Lemma 3.2]{CHR1} there exists a clique $C\in \cc$, and factorizations $g=g_1s_1w$, $h=w\iv s_2h_1$  where $s_1, s_2\in C, \lambda(s_1)=\lambda(s_2)=\lambda(s_1s_2)=|C|$ and $q=|C|+2\lambda(w)$. Let $s=s_1s_2$. Then the representation of $k$  as $g_1sh_1$ is a factorization.

Define the map $\rc\colon G^2\to G\cdot G_C^3$ by $$\rc(g,k)=(g_1,g_1s_1, g_1s_1s_2).$$ Let us prove $(rc_1)-(rc_4)$. We keep the above notation.

To prove $(rc_1)$ let us fix $k\in G$ and vary $g\in G$, $L(g)\le r$. Let $\rc(g,k)=(g_1,g_1s_1,g_1s)$. Then $\lambda(g_1)\le\lambda(g)\le r$.  Since $k=g_1sh_1$ is a factorization, the number of possibilities for $g_1$ and $g_1s$ can be bounded from above by some polynomial $P_1(r)$ by \cite[Lemma 3.1]{CHR1}. Therefore the number of possibilities for the pair of points $(g_1,g_1s)$ is bounded by, say, $P_1(r)^2$.

To prove $(rc_2)$, fix $g$, $L(g)\le r$,  and vary $k$. Let $\rc(g,k)=(g_1,g_1s_1, g_1s)$. Then since $g=g_1s_1w$ is a factorization, and $\lambda(g_1)\le r$, the number of possibilities for the pair $(g_1,g_1s)$ can be again bounded from above by $P_1(r)^2$.

To prove $(rc_3)$, fix $h=g\iv k$ and vary $g$, $L(g)\le r$. Let $\rc(g,gh)=(g_1,g_1s_1, g_1s)$. Then $$g\iv \rc(g,gh)=(g\iv g_1, g\iv g_1s_1, g\iv g_1s)=(w\iv s_1\iv, w\iv, w\iv s_2)$$ where $h=w\iv s_2h_1$ is a factorization. Again, since $\lambda(w)\le r$, the number of choices for $(w\iv, w\iv s_2)$ is bounded from above by $P_1^2(r)$. This gives $(rc_3)$.

The property $(rc_4)$ follows from the fact that if $\rc(g,k)=(g_1, g_1s_1, g_1s)$, then $g_1sw$ (for some $w$) is a factorization of $g$ (hence every syllable of $g_1s_1$ is a syllable of $g$), $g_1sh_1$ is a factorization of $k$ and $w\iv s_2h_1$ is a factorization of $g\iv k$ where $s_2=s_1\iv s$.
\endproof

Theorems \ref{t:c}, \ref{t:ab} and Remark \ref{r:d} imply the following corollary.

\begin{cor} Any graph product of groups with centroid property has the centroid property.
\end{cor}

\section{Groups without RD and open problems}\label{s:4}

Here we shall give a simple example of a group which does not contain amenable subgroups of superpolynomial growth and does not have property RD.

First let us deduce a simple ``non-amenability-like" property from RD. Suppose that a countable group $G$ with length function $L$ has property RD for some polynomial $P(r)$. Let $r>0$ and let $S, X$ be any two finite subsets of $G$, all elements of $S$ have length at most $r$. Let $\phi$ be the indicator function of $S$, let $\psi$ be the indicator function of $X$. For every $g\in SX$ let $n_g$ be the number of decompositions $g=sx$, $s\in S, x\in X$. Then
it is easy to compute that $||\phi*\psi||^2=\sum_{g\in SX} n_g^2$, $||\phi||=\sqrt{|S|}, ||\psi||=\sqrt{|X|}$. Thus
$$\sum_{g\in SX} n_g^2\le P(r)|S| |X|.$$  By (\ref{5}), we have
$$P(r)|S| |X|\ge \frac{\left(\sum_{g\in SX}n_g\right)^2}{|SX|}=\frac{|S|^2|X|^2}{|SX|}.$$
Hence we deduce

\begin{proposition}\label{l5} For every countable group $G$ satisfying property RD with respect to a length function $L$ and polynomial $P(r)$, every two finite subsets $S, X$ from $G$ such that every element of $S$ has length at most $r$ we have
$$|SX|\ge \frac{|S| |X|}{P(r)}.$$
\end{proposition}

\begin{remark}\label{r:101} It is worth noting that if we do not have the restriction that $\phi$ and $\psi$ take only values 0 and 1, we will not get a stronger inequality. Indeed, let $\phi$ and $\psi$ be arbitraty functions $G\to \R_+$ with finite supports $S, X$. Then for every $\alpha>0$ let $S_\alpha=\phi\iv(\alpha)$,  $X_\alpha=\psi\iv(\alpha)$. Then we have $||\phi||^2=\sum\limits_\alpha \alpha^2 |S_\alpha|$, $||\psi||^2=\sum\limits_\beta \beta^2|X_\beta|$. Also $$||\phi*\psi||^2=\sum\limits_{z\in SX}\left(\sum\limits_{\alpha,beta}\alpha\beta n_{\alpha,\beta}(z)\right)^2$$ where $n_{\alpha,\beta}(z)$ is the number of decompositions $z=sx, s\in S_\alpha, x\in X_\beta$.  Applying (\ref{5}), we obtain the following inequality

$$\begin{array}{l}||\phi*\psi||^2\ge \frac{\left(\sum\limits_{z\in SX}\sum\limits_{\alpha,\beta}\alpha\beta n_{\alpha,\beta}(z)\right)^2}{|SX|}=\frac{\left(\sum\limits_{\alpha,\beta}\alpha\beta|S_\alpha||X_\beta|\right)^2}{|SX|}\\ \\ =
\frac{\left(\sum\limits_\alpha\alpha|S_\alpha|\right)^2\left(\sum\limits_\beta\beta|X_\beta|\right)^2}{|SX|}.\end{array}$$ Thus property RD implies
the following inequality

\begin{equation}|SX|\ge \frac{1}{P(r)}\frac{\left(\sum\limits_\alpha\alpha|S_\alpha|\right)^2}{\sum\limits_\alpha \alpha^2 |S_\alpha|}\frac{\left(\sum\limits_\beta\beta|X_\beta|\right)^2}{\sum\limits_\beta \beta^2|X_\beta|}\label{e:101}\end{equation} for some polynomial $P$ where $r$ is the maximal length of elements in $S$. But the quotient
$\frac{\left(\sum\limits_\alpha\alpha m_\alpha\right)^2}{\sum\limits_\alpha \alpha^2 m_\alpha}$ never exceed $m=\sum\limits_\alpha m_\alpha$ if $m_\alpha\ge 0$. Indeed it is enough to write $$\left(\sum_\alpha \alpha m_\alpha\right)^2$$ as $$\left(\sum\limits_\alpha (\alpha\sqrt{m_\alpha})\sqrt{m_\alpha}\right)^2$$ and apply the Cauchy-Schwarz inequality. Thus
the right hand side of (\ref{e:101}) does not exceed $\frac {|S||X|}{P(r)}$, and
so (\ref{e:101}) follows from the inequality in Proposition \ref{l5}.
\end{remark}

Clearly, Proposition \ref{l5} shows that groups with property RD have strong expansion property. We shall call the property from Proposition \ref{l5} the \emph{property of Rapid Expansion}.
Proposition \ref{l5} immediately implies the well-known fact \cite{J} that a group with property RD cannot have an amenable subgroup of superpolynomial growth (with respect to the length function $L$ of the whole group) \cite{J}. Indeed, suppose that a group $G$ with length function $L$ has property RD and contains an amenable subgroup $H$ whose growth (with respect to $L$)  is superpolynomial. Let $r$ be such that the set $S_r$ of elements of $H$  of length $\le r$ has more than $2P(r)$ elements. Since $H$ is amenable, there exists a (F{\o}lner) set $X$ such that $|SX|<2|X|$. Then $|SX|<\frac{|S| |X|}{P(r)}$ which contradicts Proposition \ref{l5}.

In \cite[Section 3]{V1}, Valette defined another non-amenability-like consequence of the property RD. He considered radial functions (i.e., functions that are constant on spheres around the identity) instead of the indicator functions of finite sets. The property deduced in \cite{V1} is more related to the Kesten definition of amenability in terms of the spectral radius, and the Rapid Expansion property is related to the F{\o}lner definition. Of course these two approaches to amenability are close.

Using Proposition \ref{l5}, it is not difficult to construct a countable group with some length function $L$ which does not have RD and does not have amenable subgroups with superpolynomial growth. Indeed, let $G$ be the free product of all free Abelian groups $\Z^n$ of finite ranks $n\ge 1$.
Suppose that the free factor $\Z^n$ is generated by $a_{1,n},\ldots,a_{n,n}$. Let us assign to each $a_{i,n}$ the weight $n\le w_{i,n}\le p(n)$ for some increasing function $p\colon \N\to \N$. For every element $g\in G$ let $L_p(g)$ be the smallest weight of a word in $a_{i,n}$ representing $g$. Then $L_p$ is a length function on $G$. Note that $\Z^n$ has growth function $\ge C_nr^n$ for some constant $C_n$ with respect to the generating set $a_{1,n},\ldots, a_{n,n}$. Since the weight of every letter $a_{i,n}$ is at most $p(n)$,  the weight of every word in $a_{1,n},\ldots, a_{n,n}$ of length $\le r$ is at most $p(n)r$. Therefore the growth function of each $\Z^n$ with respect to the length function $L$ is at least $\left(\frac{r}{p(n)}\right)^n$. Now take any polynomial $P(r)$ of degree $n$. Then the growth function of $\Z^{n+1}$ with respect to the length function $L$ is greater than $2P(r)$ for all sufficiently large $r$. Let $S_r$ be the ball (with respect to $L_p$) or radius $r$ in $\Z^{n+1}$, $r\gg 1$.  Since $\Z^{n+1}$ is amenable,
there exists a finite set $X\subset \Z^{n+1}$ such that $|SX|\le 2|X|$. Then
$$|SX|\le 2|X|<\frac{|S| |X|}{P(r)},$$ hence $G$ with the length function $L$ does not have property RD by Proposition \ref{l5}. By Kurosh's theorem every subgroup of $G$ either contains a free non-Abelian subgroup or is conjugated to a subgroup of one of the $\Z^n$, and hence has polynomial growth (since $L_p(a_{n,i})\ge 1$).

There are also finitely generated groups without property RD and without amenable subgroups of superpolynomial growth. One such group was constructed by Denis Osin in 2012 (unpublished).

Alternatively using \cite{O99, Ol14}, one can construct a 2-generated group without property RD where every amenable subgroup is cyclic and is quasi-isometrically embedded (hence has linear growth). For this, let $G$ be the free group of countable rank freely generated by $x_i,i\in \N$. Define the length function $L$ by setting $L(x_i)=L(x_i\iv)=\log_2 i+1, i\in \N$, and for every reduced word $w$, define $L(w)$ as the sum of lengths of its letters.
It is easy to check that $L$ is a length function, and, moreover, for every $n\ge 1$ the number of elements in $G$ of length at most $n$ is at most $3^n$. The group $G$ does not have property RD by Proposition \ref{l5}.

Note that every cyclic subgroup of $G$ is obviously undistorted, hence has linear growth with respect to the length function of $G$. Let $F_2=\la a,b\ra$ be the free group of rank 2.  The group $H$ is constructed as a factor-group of the free group $F_2*G$ by the normal subgroup generated by elements $x_iw_i\iv, i\ge 1$ where $w_i$ belongs to some set of positive seventh power-free words in $a,b$ satisfying the small cancelation condition $C'(\frac{1}{100})$ such that $\log_2 i+1\le |w_i|\le d(\log_2 i+1)$ for some constant $d>1$. This set of words can be found in \cite{O99,Ol14}. Then $H$ is generated by the images of $a,b$, the natural map $x_i\mapsto w_i$ embeds $G$ into $H$, and the length function of $G$ induced by the word length function of $H$ is equivalent to $L$ (see \cite{O99}). Moreover by \cite[Theorem 1.3]{Ol14} every non-cyclic subgroup of $H$ contains a free non-Abelian subgroup and every cyclic subgroup of $H$ is undistorted by \cite[Theorem 1.4]{Ol14}.

We finish the paper with a few more open problems.

\begin{prob} Is there a finitely presented group without property RD and without amenable subgroups of superpolynomial growth?
\end{prob}

It is quite possible that one of the versions of the Higman embedding theorem (see \cite{Sap1,Sap2}) will give such a group.

\begin{prob} Is the property of Rapid Expansion equivalent to the property RD?
\end{prob}

Remark \ref{r:101} can serve as an evidence that the answer may be ``yes". In view of this remark, one can also formulate a weaker question.

\begin{prob} Is it true that property RD is equivalent to the inequality
$$||\phi*\psi||^2\le P(r)||\phi||^2||\psi||^2$$
(as in Definition \ref{d:rd})
for indicator functions $\phi$, $\psi$ of finite sets?
\end{prob}

Finally, we mention a problem first formulated by Nica \cite{N}.

\begin{prob}\label{p:nica} Find an infinite finitely generated torsion group with property RD.
\end{prob}

A formally stronger problem is to find an infinite finitely generated torsion group with the centroid property. A formally weaker problem is to find an infinite finitely generated torsion group satisfying the Rapid Expansion property.

One approach in dealing with Problem \ref{p:nica} is to consider a torsion lacunary hyperbolic group given by a presentation satisfying certain small cancelation condition considered, for example, in \cite{OOS}, and mimic the proof of Corollary \ref{c:os}.
%

\addtocontents{toc}{\contentsline {section}{\numberline { }
References \hbox {}}{\pageref{bibbb}}}

\bigskip
\begin{minipage}[t]{4.6 in}
\noindent Mark V. Sapir\\ Department of Mathematics, \\
Vanderbilt University\\
m.sapir@vanderbilt.edu
\end{minipage}

\end{document}